\documentclass[11pt,letter]{amsart}
\usepackage{amsmath, amssymb,amsthm}
\usepackage[all]{xy}
\usepackage{enumitem}

\usepackage{amsmath, amssymb,amsthm}
\usepackage{mathrsfs}

\usepackage{etoolbox}

\makeatletter
\let\old@tocline\@tocline
\let\section@tocline\@tocline
\newcommand{\subsection@dotsep}{3}
\patchcmd{\@tocline}
  {\hfil}
  {\nobreak
     \leaders\hbox{$\m@th
       \mkern \subsection@dotsep mu\hbox{.}\mkern \subsection@dotsep mu$}\hfill
     \nobreak}{}{}

\let\subsection@tocline\@tocline
\let\@tocline\old@tocline

\let\old@l@subsection\l@subsection

\def\@tocwriteb#1#2#3{%
  \begingroup
    \@xp\def\csname #2@tocline\endcsname##1##2##3##4##5##6{%
      \ifnum##1>\c@tocdepth
      \else \sbox\z@{##5\let\indentlabel\@tochangmeasure##6}\fi}%
    \csname l@#2\endcsname{#1{\csname#2name\endcsname}{\@secnumber}{}}%
  \endgroup
  \addcontentsline{toc}{#2}%
    {\protect#1{\csname#2name\endcsname}{\@secnumber}{#3}}}%

\newlength{\@tocsectionindent} 
\newlength{\@tocsubsectionindent}
\newlength{\@tocsectionnumwidth} 
\newlength{\@tocsubsectionnumwidth} 

\newlength{\@tocsectionvskip} 
\newlength{\@tocsubsectionvskip} 

\newcommand{\@tocsectionformat}{\bfseries}
\newcommand{\@tocsubsectionformat}{\mdseries}

\renewcommand{\l@section}{\section@tocline{1}{\@tocsectionvskip}{\@tocsectionindent}{}{\@tocsectionformat}}
\renewcommand{\l@subsection}{\subsection@tocline{2}{\@tocsubsectionvskip}{\@tocsubsectionindent}{}{\@tocsubsectionformat}}

\makeatother

\usepackage[hidelinks = true]{hyperref}

\newtheorem{thm}{Theorem}[section]
\newtheorem{prop}[thm]{Proposition}
\newtheorem{lemma}[thm]{Lemma}
\newtheorem{cor}[thm]{Corollary}

\newtheorem*{thm*}{Theorem}
\newtheorem{theoremletter}{Theorem}

\theoremstyle{definition}
\newtheorem{defn}[thm]{Definition}
\newtheorem{rem}[thm]{Remark}

\newcommand{\A}{\mathbf{A}}
\newcommand{\C}{\mathbf{C}}
\newcommand{\F}{\mathbf{F}}
\newcommand{\G}{\mathbf{G}}
\newcommand{\Q}{\mathbf{Q}}
\newcommand{\R}{\mathbf{R}}
\newcommand{\Z}{\mathbf{Z}}

\newcommand{\rN}{\mathrm{N}}

\newcommand{\cL}{\mathcal{L}}

\newcommand{\cO}{\mathcal{O}}
\newcommand{\cP}{\mathcal{P}}

\newcommand{\fg}{\mathfrak{g}}

\renewcommand{\fg}{\mathfrak{g}}

\newcommand{\wK}{\widetilde{K}}
\newcommand{\wG}{\widetilde{G}}
\newcommand{\Cl}{\mathscr{C}\!\ell}
\newcommand{\sw}{\mathsf{w}}

\newcommand{\lowerM}{\lower0.4ex\hbox{$M$}}

\DeclareMathOperator{\Aut}{Aut}
\DeclareMathOperator{\diag}{diag}
\DeclareMathOperator{\GL}{GL}
\DeclareMathOperator{\Gal}{Gal}
\DeclareMathOperator{\End}{End}
\DeclareMathOperator{\Lie}{Lie}
\DeclareMathOperator{\Ind}{Ind}

\DeclareMathOperator{\Hom}{Hom}
\DeclareMathOperator{\SL}{SL}
\DeclareMathOperator{\St}{St}
\DeclareMathOperator{\Spec}{Spec}
\DeclareMathOperator{\PGL}{PGL}
\DeclareMathOperator{\PSL}{PSL}
\DeclareMathOperator{\Tr}{Tr}
\DeclareMathOperator{\vol}{vol}
\DeclareMathOperator{\rH}{H}
\DeclareMathOperator{\rL}{L}
\DeclareMathOperator{\rU}{U}
\DeclareMathOperator{\MT}{MT}

\title[On the $\ell$-primary part of Abelian $3$-folds of Picard type]
{Uniform irreducibility of Galois action \\ on the $\ell$-primary part of Abelian $3$-folds of Picard type}
\author{Mladen Dimitrov \and Dinakar Ramakrishnan}

\address{University of Lille, Laboratoire Paul Painlev\'e, CNRS--UMR 8524,  59000 Lille, France}
\email{mladen.dimitrov@univ-lille.fr}

\address{Caltech, Mathematics 253-37,  Pasadena, CA 91125, USA}
\email{dinakar@caltech.edu}

\thanks{The first author acknowledges  support from  the I-SITE ULNE grant ANR-16-IDEX-0004. The second author thanks the Simons Foundation for support from the grant 523557.}

\newcommand{\dedicationpage}{
  \begin{center}
  \usefont{\encodingdefault}{pzc}{m}{n}
Dedicated to the memory of Yuri Manin
  \end{center}
}

\begin{document}

\begin{abstract}
Half a century ago Manin showed that given a  number field $k$ and a rational prime $\ell$,  there exists a uniform bound for the order of cyclic $\ell$-power isogenies between two non-CM elliptic curves over $k$. We generalize this to certain $2$-dimensional families of  abelian $3$-folds  with multiplication by an imaginary quadratic field.
\end{abstract}

\maketitle

\dedicationpage

 { \footnotesize
\tableofcontents

}

\section*{Introduction}

Given a prime number  $\ell$ and a number field $k$,  Manin showed  in \cite{manin}  that there exists an integer  $r=r(\ell,k)$ such that for any non-CM elliptic curve $E$ over $k$, $E[\ell^r]\simeq (\Z/\ell^r\Z)^2$ does not contain a $k$-rational line, or equivalently that the image of the reduction modulo $\ell^r$ of its $\ell$-adic Galois representation
\[\Gal_k=\Gal(\bar{k}/k)\longrightarrow \Aut_{\Z/\ell^r\Z}\left(E[\ell^r]\right)\simeq \GL(2,\Z/\ell^r\Z)\]
is not contained in a Borel subgroup. Manin's original proof can be greatly simplified using Faltings' proof of Mordell's conjecture, which came later.
In a series of papers Cadoret and Tamagawa established a definitive result regarding the uniform boundedness of the $\ell$-primary torsion for $1$-dimensional families of abelian varieties. In this paper we prove an analogous statement for certain $2$-dimensional families of  abelian $3$-folds which we believe to be the first  result over a genuine two-dimensional base.

Henceforth we fix an imaginary quadratic field $M$ of odd fundamental discriminant $-D\ne -3$ and denote by  $\cO_M$  its ring of integers. An abelian $3$-fold of {\it Picard type} over a field $k$ containing $M$ will always  stand for a
principally polarized abelian variety over $k$ of dimension $3$   having multiplication by $\cO_M$ defined over $k$.   Its $\ell$-adic Tate module $T_\ell A$ is free of rank $3$ over $\cO :=\Z_\ell\otimes\cO_M$   endowed with a continuous  $\cO$-linear action of $\Gal_k$. By a line (resp. plane) in $T_\ell A$, we   would  mean a $\cO$-submodule of rank $1$ (resp. $2$) which is a direct factor. More generally, given a positive integer $r$, a line (resp. plane) in $A[\ell^r]$ will always be assumed be the image, under the natural reduction map, of a line (resp. plane) in  $T_\ell A$. Finally, by a full flag we would mean a tuple of a line  sitting as direct factor in a plane.
Lines (resp. planes) will be called $k$-rational if they are stable by $\Gal_k$ (but not necessarily point-wise fixed).

Our first main result addresses  the semi-stable case.

\begin{theoremletter} \label{theoremA}
Given a number field $k$, a prime number $\ell$ inert in $M$ and a finite set $S$  of places of $M$, there exists an integer $r=r(\ell,k,S)$ such that for any non-CM abelian $3$-fold  $A$ over $k$ of  Picard type which is semi-stable outside $S$,    $A[\ell^r]$ does not contain a  full $k$-rational flag.
\end{theoremletter}

As in the case of elliptic curves, the conclusion of Theorem~\ref{theoremA} asks  the image of the attached Galois representation
\[\Gal_k\longrightarrow \Aut_{\cO/\ell^r\cO}\left(A[\ell^r]\right)\simeq \mathrm{GU}(3,\Z/\ell^r\Z)\]
not to be contained in a Borel subgroup. Also, as in the case of elliptic curves, it is necessary to cast aside the  CM abelian varieties, as their  $\ell$-adic representations are potentially  reducible.

We next show how one can relax the semi-stability assumption by adding a tiny bit of level structure at $D$.
 Given a prime $v$ of $M$ above some $p\mid D$, the projective $\Gal_k$-action on the
 $\F_p$-vector space $A[v]$ yields a homomorphism $\widetilde\rho_{A,p}: \Gal_k \to \PGL(2,\F_p)$ (see \eqref{eq:star}).
 Taking quotient by the unique index two subgroup $\PSL(2,\F_p)$ of $\PGL(2,\F_p)$ yields a canonical homomorphism
$\varepsilon_{A,p}:  \Gal_k \to\{ \pm 1 \}$ and we let $\varepsilon_{A,D}=\prod_{p \mid  D} \varepsilon_{A,p}:  \Gal_k \to\{ \pm 1 \}$.

\begin{theoremletter} \label{theoremB}
Given a number field $k$ containing $M$ and a prime number $\ell$  inert in $M$, there exists an integer $r=r(\ell,k)$ such that for any non-CM  abelian $3$-fold  $A$ over $k$ of  Picard type  and such that $\varepsilon_{A,D}$ is trivial,   $A[\ell^r]$ does not contain a full $k$-rational flag.
\end{theoremletter}

Theorem~\ref{theoremB} is the main result of this paper and implies Theorem~\ref{theoremA} as follows.
Let $k'$ be the compositum of the (finitely many) quadratic extensions of $k$ which are unramified outside
$S$ and the primes dividing $D$. Given any abelian $3$-fold  $A$ as in Theorem~\ref{theoremA}, we claim that $\varepsilon_{A,D}(\Gal_{k'})$ is trivial. Indeed, by a theorem of Grothendieck \cite[Prop.~3.5]{grothendieck-monodromie}  the semi-stability of $A$ at $v\notin S$, $v\nmid D$ implies that the inertia subgroup of $\Gal_{k}$ at $v$ acts unipotently on the $D$-adic Tate module of $A$, in particular its image by $\varepsilon_{A,D}$ is pro-$D$ hence trivial (as $D$ is odd).
Therefore the base change of $A$ to $k'$  satisfies the additional assumption in Theorem~\ref{theoremB}, implying that  Theorem~\ref{theoremA}  holds with $r(\ell,k')$ from Theorem~\ref{theoremB}.

For an individual abelian variety $A$, the  conclusion of Theorem~\ref{theoremB} is a consequence of the
Mumford--Tate conjecture which is known for abelian $3$-folds (see \S\ref{MT-groups}), so the important feature of the result is its uniformity.  As  abelian $3$-folds of Picard type are parametrized by Shimura surfaces of Picard type, a natural way to proceed would be to show that the $k$-rational points are not Zariski dense in any of their connected components $Y_\Gamma$.
Let us for the moment consider the simpler situation from  our earlier paper \cite{dimitrov-ramakrishnan} where the congruence subgroups
$\Gamma$   were neat.  Our method there had two principal steps. The first step involved showing the existence of three linearly independent global holomorphic $1$-forms on the toroidal compactification $X_\Gamma$. By a theorem of Faltings concerning the associated Albanese variety this 
 implies  that the $k$-rational points on  $X_\Gamma$ are contained in a divisor $Z$, as predicted by a conjecture of Bombieri and Lang as $X_\Gamma$ turns out to be of general type. The second step consisted in applying a result of  Nadel requiring $\Gamma$  to be neat and the  canonical divisor to be big (in his sense) to deduce that any curve $C$ of genus $\leqslant 1$ contained in $X_\Gamma$ is in fact contained in the complement of $Y_\Gamma$. Consequently, every curve in $Z$ meeting the open surface $Y_\Gamma$ must be of genus $\geqslant 2$ thus, by Faltings' proof of Mordell's conjecture for curves, $Y_\Gamma(k)$ is finite for any number field $k$. 

Let us now say a few words about the techniques involved in the  proof of Theorem~\ref{theoremB}.
 As we are  led to consider congruence subgroups of Iwahori type  $\Gamma_0(\ell^r)$, which are {\it never neat} as they have torsion, both steps mentioned above encounter difficulty and we have to resort to new methods. We produce  irregularity 
 by constructing  explicit endoscopic automorphic forms in certain  non-generic representations  $\pi$  on the unitary group in $3$ variables.  It is here  that  the  index $2$ projective Galois image condition at $D$, suggested to us by Gross,  is essential, as otherwise all of the Picard modular surfaces involved would have trivial Albanese and our approach would not apply  for global reasons. 
 Making this strategy actually work yet requires  to address  some delicate representation theoretic  questions  to which a significant part of the paper is devoted and on which we will elaborate now.

By Rogawski's theory $\pi$ is an element of an endoscopic Arthur packet  parametrized by an anti-cyclotomic (more precisely, conjugate-symplectic) Hecke character $\lambda$ of $M$.  Theorem~\ref{theoremB} imposes conditions implying that $\lambda$ must  differ from Gross' minimally ramified `canonical' characters by a finite order  character  of $M^1$
only ramified at $\ell$.  The local  Arthur packet at $\ell$ contains two representations, a supercuspidal $\pi_{c,\ell}$ and a non-tempered $\pi_{n,\ell}$, both non-generic. The difficulty of finding  $\Gamma_0(\ell^r)$-invariants in $\pi_{c,\ell}$  forces us to work with the global  $\pi_{n}$, which is automorphic if, and only if, the global root  number $W(\lambda^3)$ equals $+1$. 
For $D\equiv 3\pmod{8}$   Gross' canonical characters work and  a computation of  matrix coefficients performed in \S\ref{iwahori-inv} shows that the resulting $\pi_{n,\ell}$  has invariants even by the hyperspecial maximal compact subgroup.  

When   $D\equiv 7 \pmod{8}$  the canonical characters yield the wrong sign, leading us to consider 
$\lambda$'s which are tamely ramified   at $\ell$ to switch the sign. It remains however to  show that the non-tempered representations $\pi_{n,\ell}$ attached to such $\lambda$'s admit  $\Gamma_0(\ell^r)$-invariants for some $r$, which is significantly harder than   $\Gamma_1(\ell^r)$-invariants. To that end, in \S\ref{intertwining} we devise a more involved argument relying on Jacquet modules and intertwining operators, and which requires precise computations of exponential sums beyond the reach of known estimates. It would be worthwhile exploring  questions regarding levels of non-tempered representations in greater generality.

Once the irregularity of $X_\Gamma$ has been shown to be at least $3$ and the Bombieri--Lang conjecture established, one has to 
deal with the possible curves $C$ of genus $\leqslant 1$ contained in $Y_\Gamma$.  Using the key Lemma~\ref{lem:reflexions} on the lack of complex reflections,  our Picard modular surfaces only admit a finite number of isolated singularities hence, after removing a finite number of points,  $C$ is  endowed with an abelian family (see \S\ref{sec:abelian-family}).  This allows us apply the results of Cadoret and Tamagawa regarding the uniform boundedness of the Galois action on the Tate module of  such $1$-dimensional  families.
Finally, each of the finitely many non-CM  $k$-rational points are dealt with using  the Mumford--Tate conjecture for abelian $3$-folds of Picard type recalled in  \S\ref{MT-groups}, and which is not available in higher dimension.

As our Picard modular surfaces $X_\Gamma$ have irregularity  $\geqslant 3$, the Enriques--Kodaira classification implies that they are either ruled of genus $q$, or elliptic, or else they are of general type. In the last case, which according to Holzapfel \cite[\S 5.4]{holzapfel-book} occurs for all odd $D\notin\{3, 7,11,19, 23, 31, 39, 47, 71\}$,  we show that the   Bombieri--Lang Conjecture holds, {\it i.e.}, that  the $k$-rational points are not Zariski dense.  Investigating small values of $D$, as suggested by Mazur,  seems  even more interesting. It is established in {\it loc. cit.} that for all $D\ne 71$ in the above list the level $1$ Picard modular surfaces are rational and it would be natural to investigate the nature of their degree $2$ Gross covers that we consider. A way to shed light on this question  would be to find an explicit $2$-parameter family of abelian $3$-folds of  Picard type to which our theorem  applies.  

It might be worthwhile remarking that we could have also considered the simpler case of the moduli of principally polarized abelian {\it surfaces} $A$ over $k$ with multiplication by $\cO_M$,  which will involve ${\rm U}(1,1)$. However, as ${\rm SU}(1,1)\simeq \SL(2)$
this essentially  reduces  to the  modular curve case. On the other hand, if we consider
 principally polarized abelian  surfaces $A$ with {\it real multiplication}, then the family is parametrized by a Hilbert modular surface which has trivial ${\rm H}^1$, thus our methods, which rely on the Albanese variety, do not lead to an establishment the  Bombieri--Lang Conjecture.

\addtocontents{toc}{\setcounter{tocdepth}{2}}

\section{Levels for endoscopic non-tempered representations of $\rU(3)$}

It goes back to the work of Casselman that  admissible irreducible  representations having  non-zero Iwahori invariants are exactly those occurring as sub-quotients in  parabolic inductions of unramified characters.  Whereas  the dimension of the invariants by the depth $r$  Iwahori subgroup  in the
 full induced representation grows as  $r$  goes to infinity, this might not always be the case for all its sub-quotients, as shown by  the example of the trivial representation of $\GL(2)$, realized as a quotient of a unramified principal series representation.

Another challenging question is to determine which sub-quotient of a parabolically induced unramified character picks up the invariants by a  given maximal open compact subgroup. Whereas MacDonald's formula for zonal spherical functions yields an answer in the case of a maximal hyperspecial subgroup, the general case appears to be an open question.

In this section we fully answer those two natural questions in the case of certain non-tempered endoscopic representations of $\rU(3)$ attached  to a quadratic  extension $E/\Q_p$. It will be later applied in a global setting to $E=M_p$, where $M$ is  an imaginary quadratic field in which the prime $p$ does not split.

{\it In this section  of our paper we will adopt local notations. }

Let $E$ be a quadratic field extension of $\Q_p$, $\cO$ be its ring of integers, $\cP$ its maximal ideal and $\varpi$ a uniformizer.
We assume that $E$ is {\it not} a ramified extension of $\Q_2$ and we  fix a generator $\xi$  of its different   ideal such that $\bar\xi=-\xi$. 
  Denote by $x\mapsto \bar x$  the  non-trivial element of $\Gal(E/\Q_p)$.
We fix an additive character $\psi: \Q_p\to \C^\times$ of conductor $0$, i.e. $\ker(\psi)=\Z_p$, and we consider the additive character
$\psi^{}_E$ of $E$ defined as $\psi^{}_E(z)=\psi(\Tr_{E/\Q_p}(z))$.

   Let $G$  be the  unique  quasi-split unitary group  in $3$ variables relative to the extension  $E/\Q_p$. It can be realized as the automorphisms of  $E^3$ preserving  the hermitian pairing
\[\langle x,y\rangle = \bar x_1 y_3 +  \bar x_2 y_2 +\bar x_3 y_1.\]

The standard Borel $B$ of $G$ is a product of its torus 
\[T=\left\{ \left(\begin{smallmatrix}\bar{\alpha}  & &   \\ &\beta & \\  & & \alpha^{-1}\end{smallmatrix} \right)\Big{|} \alpha\in E^\times, \beta\in E^1 \right\}\] with its unipotent subgroup 
\[N=\left\{ [z,x]=
\left(\begin{smallmatrix}1  & -\bar{z}& \xi x -z\bar z /2 \\ &1 & z\\  & & 1\end{smallmatrix} \right)\Big{|} z\in E, x\in \Q_p  \right\}.\]

\subsection{The Bruhat-Tits tree of  $\rU(3)$}\label{tree}
As $G$ has rank $1$, its Bruhat-Tits building is a tree. We will first describe its standard apartment.
The relative roots of $G$  are obtained by decomposing the adjoint action on the  Lie algebra of the maximal $\Q_p$-split torus $T_0=\left\{\mathrm{diag}(a,1,a^{-1}) |    a\in \Q_p^\times \right\} $ of $G$. The positive elements of the associated root system $\Phi$ are $\{\zeta, 2\zeta\}$. Let $h: \G_m \to T_0  \subset  G$ be the generator of the co-character lattice $X_\ast(T_0) \simeq \Z$ such that $\langle \zeta, h\rangle=1$. Then the co-root sub-lattice is generated by $\zeta^\vee = 2h$, so that we have the standard normalization $\langle\zeta, \zeta^\vee\rangle=2$. According to  \cite[\S 1.15]{tits} the affine roots are
$\{ \pm \zeta+\Z\}\cup \{ \pm 2\zeta+\Z\}$ if $E$ is unramified, and
$\{ \pm \zeta+\tfrac{1}{2}\Z\}\cup \{ \pm 2\zeta+\Z+\tfrac{1}{2}\}$ if $E$ is ramified;
note that $\delta=0$ in {\it loc. cit.} as $E$ is not a ramified extension of $\Q_2$.
The apartment associated to $T_0$ is $\R h$ and its  walls are the vanishing sets of these (affine) roots,
hence they are given by $\tfrac{1}{2}\Z h= \Z h \cup \tfrac{1}{2}\Z h$, resp. $\tfrac{1}{4}\Z h=
\left(\tfrac{1}{2}\Z +\tfrac{1}{4}\right)h\cup \tfrac{1}{2}\Z h$,  if  $E$ is unramified, resp. ramified.
A conjugacy class of maximal compact subgroups can be represented by a wall in the standard apartment.
By definition, a wall is hyperspecial if for every $\zeta'\in \Phi$ there exists an affine root with gradient
$\zeta'$ vanishing on that wall. Since $\left(\tfrac{1}{2}\Z +\tfrac{1}{4}\right)h\cap \tfrac{1}{2}\Z h=\varnothing $ this only can happen when $E$ is unramified, in which case the hyperspecial walls are
$ \Z h \cap \tfrac{1}{2}\Z h= \Z h$. All walls are special, as elements of $\Phi$ are rational multiple of one another. 

Let us  describe the  conjugacy classes of maximal compact subgroups in $G$ in terms equivalent classes of $\cO$-lattices $\cL$  in $E^3$ modulo homothety. Recall the definition of  the dual  lattice 
\[\cL^\perp=\Hom_{\cO}(\cL, \cO)=  \{x\in E^3 |    \langle x, \cL\rangle \subset  \cO\}.\]
There are two conjugacy classes of maximal compact subgroups in $G$, those which are  stabilizers of
self-dual lattices, and those which are stabilizers of almost-self-dual lattices ({\it i.e.}  $\cL$ such that $\cL\subsetneq \cL^\perp \subsetneq\varpi^{-1}\cL$).  
We next give an explicit description of the maximal compact subgroups corresponding to the
walls of a chamber in the standard apartment.

The standard maximal compact subgroup $K^\circ=G(\cO)$ of $G$ is defined as the stabilizer of the self-dual lattice $\cL^\circ=\cO^3$. It is hyperspecial if and only if $E$ is unramified. 
The reductive quotient $\overline{G}{}^\circ$ is given by $\rU(3,\F_p)$ if $E$ is unramified, and by $\mathrm{O}(3,\F_p)$ if $E$ is ramified.

The other standard maximal compact subgroup  $K'$ of $G$, defined as the stabilizer of the  almost self-dual  lattice $\cL'=\cO\oplus \cO\oplus \cP$, is given by
 $$K'=\left(\begin{smallmatrix} \cO & \cO& \cP^{-1}\\
\cP & \cO^\times & \cO \\ \cP   & \cP  & \cO \end{smallmatrix} \right)\cap G.$$
One has
$\cL'^\perp=\cP^{-1} \oplus \cO  \oplus \cO$  and
 $K'$ acts on $\cL'/\varpi\cL'^\perp\simeq \cO/\cP$ via its middle coefficient.
The reductive quotient $\overline{G}{}'$ is isomorphic to   $ (\rU(1,1)\times \rU(1))(\F_p)$ if $E$ is unramified, and to
$\pm \mathbf{1}_2 \cdot \SL(2,\F_p)\times\{\pm 1\}$,  if $E$ is ramified.

The standard Iwahori subgroup of $G$, defined as
$I=K^\circ\cap K'=\left(\begin{smallmatrix} \cO^\times & \cO & \cO\\
\cP  & \cO^\times &\cO \\ \cP  &\cP  & \cO^\times \end{smallmatrix} \right)\cap G$, is the stabilizer 
of a chamber in  the standard apartment in the Bruhat-Tits tree of $G$:
$$\xymatrix@C=40pt{
\ar@{--}[r]  & \gamma^{-1} K' \gamma
  \ar@{-}[r] & K^\circ \ar@{-}_{I}[r] \ar@/^1pc/@{.}[rr]^{I_{2,1}}& K'
  \ar@{-}[r] & \gamma K^\circ \gamma^{-1} \ar@{--}[r] &}$$
where  $\gamma=\left(\begin{smallmatrix}\varpi^{-1}  & & \\ & 1 & \\ & & \varpi \end{smallmatrix}\right)$
and $I_{2,1}=K^\circ  \cap \gamma K^\circ \gamma^{-1}=\left(\begin{smallmatrix} \cO^\times & \cO &  \cO\\
\cP  & \cO^\times &\cO \\ \cP^2  & \cP  & \cO^\times \end{smallmatrix} \right)\cap G$. 

One has  $G\supset K'\supset I \supset I_{2,1}\supset I_2$, where $I_r=\left(\begin{smallmatrix}\cO^\times &\cO &\cO\\\cP^r &\cO^\times &\cO\\ \cP^r & \cP^r &\cO^\times\end{smallmatrix} \right)\cap G$.
Finally, we let $K_1\subset K^\circ$ be the principal congruence subgroup of matrices  $\equiv \mathbf{1}_3\pmod{\cP}$  and  
$K_T=\left(\begin{smallmatrix}\cO^\times &\cP &\cP\\\cP &\cO^\times &\cP\\ \cP & \cP&\cO^\times \end{smallmatrix} \right)\cap G$.

\subsection{Review  of $L$-parameters and $A$-packets} \label{packets}

For any integer $n\geqslant 1$ there are exactly two (up to isomorphism) $n$-dimensional hermitian spaces over $E$, depending on the image of the discriminant in $\Q_p^\times/\rN_{E/\Q_p}(E^\times)$, and the corresponding unitary groups $\rU(n)$ are isomorphic if and only if $n$ is odd. When $n=2$, by analogy with the Archimedean case, we will denote by $\rU(1,1)$ the quasi-split form and by $\rU(2)$ the compact one.

The $L$-group   of $\rU(n)$ is given by $\GL(n,\C)\rtimes W_{\Q_p}$ with the Weil group $W_{\Q_p}$ acting on $\GL(n,\C)$ through its
quotient $\Gal(E/\Q_p)$ whose non-trivial element sends  $g$ to $w_n{}^{t}g^{-1} w_n^{-1}$, where   $w_n$   denotes the  anti-diagonal matrix $(1,-1,1,\dots,(-1)^{n-1})$. By definition, an $L$-parameter for  the  quasi-split  $\rU(n)$  is  a homomorphism
$W_{\Q_p}\times \SL(2,\C) \longrightarrow \GL(n,\C)\rtimes W_{\Q_p}$,  but as one knows (see \cite[\S 3]{GGP1})  it is equivalent to
ask for its restriction
\[\phi: W_E\times \SL(2,\C) \longrightarrow \GL(n,\C), \]
to be {\it conjugate-orthogonal} if  $n$ is odd and {\it conjugate-symplectic} if $n$ is even.
Recall that $\phi$  is conjugate-self-dual  if $\overline\phi\simeq \phi^\vee$, or equivalently, if the induced representation $\Ind_{W_E}^{W_{\Q_p}}(\phi)$ is self-dual. Furthermore,  $\phi$ is
conjugate-orthogonal, resp. conjugate-symplectic, if it  preserves a non-degenerate symmetric, resp. skew-symmetric, bilinear form. Note that while Schur's Lemma implies that any irreducible self-dual  (or conjugate-self-dual) parameter has a well defined sign, this need not be always the case for reducible parameters.

For $n=1$, a character of  $E^\times$ is conjugate-orthogonal (resp. conjugate-symplectic) if its
 restriction to $\Q_p^\times$ is   trivial (resp. is the quadratic character attached to $E/\Q_p$). For $n\in \Z_{\geqslant 0}$, the $n$-th symmetric power of the standard $2$-dimensional representation $\St$ of SL$(2,\C)$, with $W_E$ acting trivially, is  conjugate-symplectic if $n$ is odd and conjugate-orthogonal if $n$ is even.

The base change $\nu^{}_E(z)=\nu(z/\overline z)$ from $\rU(1)$ to $\GL(1)/E$ of a character $\nu$  of $E^1$ is conjugate-orthogonal and conversely any 
conjugate-orthogonal character of $E^\times$ is  obtained in that way. 
For $\lambda$   a conjugate-symplectic character of $E^\times$, the  conjugate-orthogonal representation
\[  (\lambda \otimes \St) \oplus  \nu^{}_E  : W_E\times \SL(2,\C) \longrightarrow \GL(3,\C) \]
would be of key relevance to us. It yields an  $L$-parameter $\phi_{\lambda, \nu}$ of $G$, coming from an $L$-parameter of the (unique) cuspidal endoscopic subgroup $H= \rU(1,1) \times \rU(1)$ of $G$. The cardinality of the corresponding $L$-packet  $\Pi_L(\phi_{\lambda, \nu})$ is given by  the order of the centralizer of $\phi_{\lambda, \nu}$ (modulo  center) which turns out to be $2$. More precisely, $\Pi_L(\phi_{\lambda, \nu})$ contains two  discrete series representations $\pi_2$ and $\pi_c$ of $\rU(3)$, exactly one of them,  namely $\pi_c$, being supercuspidal (see  \cite[Chap.~12.2]{rogawski-U3} where this $L$-packet is denoted $\Pi_L(\St_H(\xi))$). There is another endoscopic $L$-packet for $G$ consisting of a single 
non-tempered representation $\pi_n$ whose  the $L$-parameter  is given by
\[\lambda |\cdot |_E^{1/2}\oplus \lambda |\cdot|_E^{-1/2}\oplus \nu^{}_E: W_E\times \SL(2,\C) \longrightarrow \GL(3,\C). \]

 Rogawski's theory \cite{rogawski-U3,rogawski-A-packets} describes the automorphic representations contributing to the $\rH^1$ of  Shimura surfaces of Picard type  in terms global Arthur packets (see \cite[\S 3.1]{dimitrov-ramakrishnan} for a summary). The corresponding local Arthur packet  at $p$ has $2$ elements
 $\Pi(\lambda,\nu)=\{\pi_n,\pi_c\}$ (see \cite[\S 12.3.3]{rogawski-U3}, where $\pi_c$ is denoted $\pi^s$), and the restriction to $W_E$ of its  $A$-parameter  is given by
\[
(\lambda\otimes \mathbf{1} \otimes \St) \oplus  \nu^{}_E  : 
W_E\times \SL(2,\C)\times \SL(2,\C) \longrightarrow \GL(3,\C), 
\]
while the $A$-parameter of $\pi_2$ is given by $(\lambda  \otimes \St \otimes \mathbf{1}) \oplus  \nu^{}_E$.

Crucial for us would be the description $\pi_n$ and $\pi_2$ as  the Jordan--H\"older constituents of a 
principal series representation $\pi$. Indeed, by \cite[\S 1]{rogawski-A-packets}, $\pi_n$ is the Langlands quotient 
of the  (unitarily normalized) parabolic   induction  of the character 
\begin{equation}\label{eq:mu}
\mu(\bar\alpha, \beta, \alpha^{-1})= \lambda(\bar\alpha)\nu(\beta)|\alpha|_E^{1/2}, 
\end{equation}
with $\pi_2$  the unique non-zero irreducible sub-representation. The sub and quotient are switched when $\mu$ is  replaced by 
$\mu^{\sw}(\bar\alpha, \beta, \alpha^{-1})= \lambda(\bar\alpha)\nu(\beta)|\alpha|_E^{-1/2}$,
 where $\sw=\left(\begin{smallmatrix}  & &  1 \\ &1 & \\ 1 & & \end{smallmatrix} \right)$ is the non-trivial element of the Weyl group of $G$. 
The Jacquet functor $\pi\to \pi_N$ is exact and it sends $\pi_2$ (resp. $\pi_n$) to 
$\mu\delta^{1/2}$ (resp. $\mu^{\sw}\delta^{1/2}$), where $\delta(\bar\alpha, \beta, \alpha^{-1})= |\alpha|_E^2$ is the modulus character.  The extension 
\begin{equation}\label{induction}
0\to \pi_2 \to \pi = \Ind_B^G(\mu)\xrightarrow{\mathrm{pr}} \pi_n \to 0
\end{equation} 
does not split since, by Frobenuis reciprocity, one has
\[\End_G\left(\Ind_B^G(\mu)\right)=\Hom_B\left(\Ind_B^G(\mu), \mu\delta^{1/2}\right)=
\Hom_T\left(\Ind_B^G(\mu)_N, \mu\delta^{1/2}\right)\simeq \C. \]
One knows by Rodier  \cite[Thm.~2]{rodier}  that the image of $\Ind_B^G(\mu)$ by a twisted (by a non-degenerate character of $N$) Jacquet functor, singling out generic representations,  is a line.  Since  $\pi_n$ is non-generic  (see \cite[p.174]{rogawski-U3}), the exactness of this functor implies that $\pi_2$ is generic.

As $\pi_n$ is non-tempered, the subspace $\pi_2$ consists of  $f\in\pi $ such that for all $f^\vee\in \pi ^\vee$ the matrix coefficient 
$g\mapsto \langle g\cdot f,  f^\vee \rangle$ belongs to $\rL^2(G)$. Conversely the following lemma holds. 

\begin{lemma} \label{lem:L2-crit} Let $f\in \pi $. 
If  $g\mapsto \langle g\cdot f,  f^\vee \rangle$ belongs to $\rL^2(G)$ for some $0\ne f^\vee\in \pi ^\vee$, then $f\in \pi_2$. 
\end{lemma}
\begin{proof} The dual of \eqref{induction} is given by
\[0\to \pi_n^\vee \to \pi ^\vee=\Ind_B^G(\mu^{-1}) \to \pi_2^\vee \to 0,\]
and the irreducibility of  $\pi_2$ and $\pi_n$ implies that $\pi_n^\vee=\{f^\vee\in \pi ^\vee | \langle \pi_2, f^\vee\rangle=0\}$. As  $f^\vee\ne 0$, its $G$-span contains  $\pi_n^\vee$, implying that the matrix coefficient  
$g\mapsto \langle g\cdot f,  f^\vee \rangle$ belongs to $\rL^2(G)$  for all $f^\vee\in\pi_n^\vee$. One deduces that
\[g\mapsto \langle g\cdot f,  f^\vee \rangle= \langle \mathrm{pr}(g\cdot f),  f^\vee \rangle=\langle g\cdot \mathrm{pr}(f),  f^\vee \rangle\in \rL^2(G)\]
As the irreducible $\pi_n$ is not a discrete series representation, this implies  $\mathrm{pr}(f)=0$, {\it i.e.} $f\in \pi_2$. 
   \end{proof}

We will be mostly interested in the following $A$-packets having trivial central characters:
\begin{equation}\label{A-packet-lambda}
\Pi(\lambda)=\Pi(\lambda,\lambda_{| E^1}^{-1}).
\end{equation}

\subsection{The Gross subgroup $K''$}\label{K-second}
In this  subsection, $E$ is  assumed ramified (hence $p$ is  odd),   $\cO/\cP=\F_p$ and $\cP=(\xi)$. 
As $|\PGL(2,\F_p)|=|\SL(2,\F_p)|$ all vertices in the  tree of $G$  have valence $p^3+1$.
 The reductive quotient $\overline{G}{}^\circ$ of $K^\circ$  is isomorphic to  the orthogonal group  $\mathrm{O}(3,\F_p)$ with respect to the quadratic form represented by   $\left(\begin{smallmatrix}  & & 1 \\ & 1 & \\ 1 & &  \end{smallmatrix}\right)$. The  adjoint action on matrices $\left(\begin{smallmatrix} y & x \\  z  & -y \end{smallmatrix}\right)$
preserving the determinant $-(y^2+xz)$  allows us to identify $\PGL(2,\F_p)$ and $\mathrm{SO}(3,\F_p)$   as follows:
\begin{equation}\label{eq:adjoint}
\begin{pmatrix} a & b \\  c  & d \end{pmatrix}\mapsto \frac{1}{ad-bc}
\begin{pmatrix} a^2 & -ab & -b^2/2\\ -2ac  & ad+bc &bd \\ -2c^2  &2cd  & d^2 \end{pmatrix}.
\end{equation}
One has $\mathrm{O}(3,\F_p)=\pm \mathbf{1}_3 \cdot \mathrm{SO}(3,\F_p)$. By the above description,   $\mathrm{SO}(3,\F_p)$ is generated by the~set
\[\left\{\left(\begin{smallmatrix} & & -1/2\\   & -1 &\\ -2 & & \end{smallmatrix}\right),
\left(\begin{smallmatrix} a & &\\  & 1 & \\ & & a^{-1} \end{smallmatrix}\right),
\left(\begin{smallmatrix} 1 & -b & -b^2/2\\   & 1 &b \\   &  & 1 \end{smallmatrix}\right)\Big{|}   a\in \F_p^\times, b\in\F_p\right\}.\]

\begin{defn}\label{def-double-prime}
Let $K''$  be the index $2$ subgroup of $K^\circ$ defined as the inverse image of the
subgroup of $\mathrm{O}(3,\F_p)$ generated by  $-\mathbf{1}_3$  and the image of $\PSL(2,\F_p)$. Let $I''=K''\cap K'\subset I$.
\end{defn}

As $E/\Q_p$ is a ramified quadratic extension with $p$ odd, a conjugate-symplectic character 
$\lambda$ of $E^\times$ is necessarily ramified and its restriction to $\Z_p^\times$ is given by its unique quadratic character. 
If $\lambda$ is tamely ramified, then  its restriction to $\cO^\times$ is also given by its unique quadratic character, and 
the equation $ \lambda (\xi)^2= \lambda(-\xi\bar\xi)=\lambda(-1)=(-1)^{(p-1)/2}$ shows that there are precisely two such characters.  

Interested in determining a level for an element of the $A$-packet $\Pi(\lambda)$ considered in \eqref{A-packet-lambda}, 
we are indebted to B.~Gross for generously sharing a suggestion that led to the following proposition.

\begin{prop} \label{ram-inv}  
Let $\lambda$ be a  tamely ramified conjugate-symplectic character of $E^\times$ and let 
$\pi_{n}$ be the non-tempered member of the $A$-packet $\Pi(\lambda)$. 
Then  $\dim\pi_{n}^{K''}=\dim\pi_2^{I''}=1$.
\end{prop}
\begin{proof} As $\mu_{|(T\cap I'')}=\mu^{\sw}_{|(T\cap I'')}=\mathbf{1}$, applying the  Jacquet functor to the exact sequence
of admissible  $G$-representations \eqref{induction} allows one to see that both  $\pi_2$ and $\pi_n$ have non-trivial $I''$-invariants.
Moreover, as $\big{|}B\backslash G/ I'' \big{|}=\big{|}(B\cap K'')\backslash K''/ I'' \big{|}=2$,  both $\pi_2^{I''}$ and $\pi_n^{I''}$ must be   $1$-dimensional.
By Iwasawa decomposition, the restriction of $\Ind_B^G(\mu)$  to $K''$ is given by $\Ind_{B\cap K''}^{K''}(\mu)$, hence  the
line $\Ind_B^G(\mu)^{K''}$ admits a basis  $f$ uniquely characterized by $f_{|K''}= \mathbf{1}_{K''}$. It follows that
 $\dim\pi_{n}^{K''}+\dim\pi_2^{K''}=1$ and we will show that $\pi_2^{K''}=\{0\}$.

The  line $(\Ind_B^G(\mu^{-1}))^{K''}$ admits a basis  $f^\vee$ uniquely characterized by
\[\langle v, f^\vee\rangle = \frac{1}{\sqrt{\vol(K'')}} \int\limits_{K''} v(k) dk. \]

By Lemma~\ref{lem:L2-crit} one has   $f\notin \pi_2$ if, and only if,
$(g\mapsto \langle g\cdot f,  f^\vee \rangle)\notin \rL^2(K''\backslash G/K'')$.

Recall    $\gamma =\left(\begin{smallmatrix} -\xi^{-1}& &  \\ & 1 & \\ & &  \xi  \end{smallmatrix} \right)$ and let
$\eta =\left(\begin{smallmatrix} \bar u& &  \\ & 1 & \\   & & u^{-1} \end{smallmatrix} \right)$ where
$u\in\cO^\times$ is a fixed non-square element.

As $K^\circ=K''\coprod \eta K''$,  Cartan decomposition for the special maximal compact $K^\circ$ yields:
\[ G= \coprod_{n\geqslant 0} \left(K''\gamma^n K''\right) \amalg \left(K'' \gamma^n \eta K''\right).\]

Since $\eta\cdot f=-f$ one deduces that  $\langle  \gamma^n \eta \cdot f,  f^\vee \rangle=-\langle \gamma^n \cdot f,  f^\vee \rangle$
and checking  that $f\notin \pi_2$  amounts to proving the divergence of  the numerical sequence with general term
\[\vol(K''\gamma^n K'') \big|\langle \gamma^n \cdot f,  f^\vee \rangle \big|^2=[K'':(K''\cap \gamma^n K''\gamma^{-n})]
 \Big|\int\limits_{K''} f(k \gamma^n) dk \Big|^2.\]

By Iwahori decomposition one has $[I'':(K''\cap \gamma^n K''\gamma^{-n})]=p^{2n-1}$ for all $n\geqslant 1$.
As $\big|(\mu\delta^{1/2})(\gamma)\big|= p^{3/2}$ we are led  to establish the divergence of  the sequence with general term
\[ \Phi_n=p^{5n/2} \cdot \Big|\int\limits_{K''} f(\gamma^{-n} k \gamma^n) dk \Big|.\]

In view of the  inequality $p>\sqrt{p}+1$ for $p\geqslant 3$, this will follow from the next lemma.  \end{proof}

\begin{lemma} For all $n\geqslant 1$ one  has  $\left|\int\limits_{K''{\smallsetminus} K''_{2n}} f(\gamma^{-n} k \gamma^n) dk \right|\leqslant \vol(I'') \cdot (\sqrt{p}+1)p^{-2n}$ and    $\left|\int\limits_{K''_{2n}} f(\gamma^{-n} k \gamma^n) dk \right|=\vol(I'')\cdot p^{1-2n}$.
\end{lemma}

\begin{proof}
The last row of an element $k\in K''$ is given by  $(0,0,1)\cdot k= (c_1(k),  c_2(k), c_3(k)) \in \cO^3{\smallsetminus}   \cP^3$. For $j\geqslant 0$ we let
\[K''_j=\left\{k\in K'' \big{|} c_1(k)\in \cP^j \right\} \text{ and } K^{\prime\prime\times}_j=K''_j{\smallsetminus} K''_{j+1}.\]
Note that $K''_0=K''$ and $K''_1=I''$. We use the partition $K''{\smallsetminus} K''_{2n}=K_0^{\prime\prime\times}\coprod K^{\prime\prime\times}_1\coprod \dots \coprod  K^{\prime\prime\times} _{2n-1}$   to compute the first integral and $K''_{2n}=K^{\prime\prime\times}_{2n} \coprod K''_{2n+1}$
for the second.

For $j\geqslant 1$, using Iwahori decomposition,  one finds that $[I'':K''_j]=  
[I''\cap N^{-}: K''_j\cap N^{-}]=p^{j-1}$.

Given any $k\in K^{\prime\prime\times}_j (0\leqslant j \leqslant 2n)$, using the Iwasawa decomposition $G=\gamma^{\Z} N K^\circ=\gamma^{\Z} N K'$,
one finds that $\gamma^{n-j} k \gamma^n\in N K^\circ$, hence $\vert f(\gamma^{-n} k \gamma^n) \vert\leqslant \vert \mu(\gamma^{j-2n})\vert=
p^{\frac{3}{2}j-3n}$.  Therefore
\[\left|\int\limits_{K''{\smallsetminus} K''_{2n}} f(\gamma^{-n} k \gamma^n) dk \right|\leqslant \vol(K''{\smallsetminus} I'')\cdot p^{-3n} + \vol(I'')
\sum_{j=1}^{2n-1} p^{\frac{1}{2}j-3n}(p-1)= \]
\[=\vol(I'') \left(p^{\frac{1}{2}-2n}+p^{-2n}-p^{1-3n}-p^{\frac{1}{2}-3n}+([K'':I'']-1)p^{-3n}   \right),\]
proving the desired inequality, as $[K'':I'']=p+1$ (obtained by  going to the reductive quotient).

Since $\vol(I'')\cdot p^{1-2n}=\vol(K''_{2n})$, in order to complete the proof of the lemma,  it suffices to show that
$f(\gamma^{-n} k \gamma^n)$ is constant on  $k\in K''_{2n}$. This is evident for $k\in K''_{2n+1}$, as then $k$ and
$\gamma^{-n} k \gamma^n$ both belong to $I''$ and share same determinant and lower right coefficient $c_3$, implying that
$f(\gamma^{-n} k \gamma^n)=f(k)$. Miraculously, as one can see from \eqref{eq:adjoint}, this remains true for $k\in K''_{2n}{\smallsetminus} K''_{2n+1}$ as well, {\it i.e.} even though $\gamma^{-n} k \gamma^n\in K^\circ {\smallsetminus} I$, the fact that $c_3(\gamma^{-n} k \gamma^n)=c_3(k)$ still implies that
$\gamma^{-n} k \gamma^n\in K''$.  \end{proof}

\subsection{Higher Iwahori invariants via matrix coefficients}\label{iwahori-inv}
In this  subsection,  we assume that $E/\Q_p$ is unramified and we let 
$\lambda_0$ denote  the unique quadratic unramified character of $E^\times$. By  Iwasawa decomposition, the corresponding  $\pi=\Ind_B^G(\mu^{}_0)$ has one dimensional invariants by any given maximal open compact subgroup $K$  of $G$.   The following proposition  states, depending on $K$,
 whether the $K$-invariant line belongs its sub-representation $\pi_2$ or maps non-trivially to 
 the quotient  $\pi_n$ (see \S\ref{packets} for notation).
We recall that $K^\circ$ and $K'$ are the two standard maximal compact subgroups, $K^\circ$ being the hyperspecial one, and that the standard Iwahori subgroup $I$ equals $K^\circ \cap K'$. 

\begin{prop}\label{prop:inert-inv} 
One has $\dim\pi_2^{K'}=\dim\pi_n^{K^\circ}=1$.
\end{prop}

\begin{proof}
Applying the  Jacquet functor to the exact sequence \eqref{induction} allows one to see that both $\pi_2$ and $\pi_n$ have $1$-dimensional $I$-invariants.  As by Cartan decomposition $G$ is generated by $K$ and $\gamma$, hence by $K^\circ$ and $K'$, it follows  that exactly  one amongst $\pi_2^I$ and  $\pi_n^I$  is fixed by $K^\circ$, while the other one is fixed by $K'$.  By Iwasawa decomposition, the restriction to $K$  yields an isomorphism  $\pi\xrightarrow{\sim}\Ind_{B\cap K}^K(\mu^{}_0)$, under which a  basis  $f_K$  of $\pi ^K$ is mapped to 
$\mathbf{1}_K$. Moreover, the  line $(\pi ^\vee)^K$ admits a basis  $f_K^\vee$  characterized by
\[\langle f, f_K^\vee\rangle = \frac{1}{\sqrt{\vol(K)}}\int\limits_K f(k) dk. \]

The remainder of the proof consists in computing the bi-$K$-invariant function $g\mapsto \langle g\cdot f_K,  f_K^\vee \rangle$ and checking whether it belongs or not to
$\rL^2(K\backslash G/K)$. Using Cartan decomposition $G=\coprod_{n\geqslant 0} K\gamma^n K$ this amounts to checking whether
$\rL^2(\Z_{\geqslant 0})$ contains the numerical sequence
\[\sqrt{\vol(K\gamma^n K)} \langle \gamma^n \cdot f_K,  f_K^\vee \rangle=\sqrt{[K:(K\cap \gamma^{-n}K\gamma^n)]}
\int\limits_K f_K(k \gamma^n) dk.\]
Using Iwahori decomposition for all $n\geqslant 1$ we have  $[K:(K\cap \gamma^{-n}K\gamma^n)]/[K:I]=p^{4n-3}$ (resp. $p^{4n-1}$), where  $K=K^\circ$ (resp. $K'$). Since $(\mu^{}_0\delta^{1/2})(\gamma)= -p^3$ we have $f_K\in \pi_2$ if, and only if,
\begin{align}\label{Phi-n}
(\Phi^K_n)_n\in \rL^2(\Z_{\geqslant 0}) \text{, where }\Phi^K_n=p^{5n} \cdot \int\limits_K f_K(\gamma^{-n} k \gamma^n) dk.
\end{align}
The proof of the Proposition is then completed  by the following  Lemma.
\end{proof}

\begin{lemma} The quantity  $p^{n}\cdot\Phi^{K'}_n$ is independent of $n\geqslant 1$, in particular  $(\Phi^{K'}_n)_n\in \rL^2(\Z_{\geqslant 0})$.
\end{lemma}
\begin{proof}
The last row of an element $k\in K'$ is given by  $(0,0,1)\cdot k= (p\cdot c_1(k), p\cdot c_2(k), c_3(k))$ with $c_2(k)\in \cO$ and $(c_1(k),c_3(k))\in
(\cO\times \cO) {\smallsetminus} (\cP\times \cP)$. For $j\geqslant 0$ we let
\[K'_j=\left\{k\in K' \Big{|} c_1(k)\in \cP^j\right\} \text{ and } K^{\prime\times}_j=K'_j{\smallsetminus} K'_{j+1}.\]
To compute the above integral we use the partition $K'=K_0^{\prime\times}\coprod K^{\prime\times}_1\coprod \dots \coprod  K^{\prime\times} _{2n-1}\coprod K'_{2n}$.

First, we compute the volume of $K'_j$, for $j\geqslant 1$. Using Iwahori decomposition one finds that
\[ [K':K'_j]= \frac{[K':I]}{[K'_j:I\cap K'_j]}[I:I\cap K'_j]=  \frac{[K':I]}{[K' \cap N : I\cap N ]}[I\cap N^{-}: K'_j\cap N^{-}]=c_0\cdot p^{j+2\left[\frac{j}{2}\right]}. \]

Next we observe that  by Iwasawa decomposition, for all $k\in K'_{2n}$ one has $c_2(k)\in \cP^n$, {\it i.e.} $\gamma^{-n} k \gamma^n\in N\cdot K'$, and
therefore $f_K(\gamma^{-n} k \gamma^n)=1$.

Using again Iwasawa decomposition, one checks that  for  $0\leqslant j \leqslant 2n-1$ and for every $k\in K^{\prime\times}_{j}$ one has
$p^{n-j} c_2(k)\in \cO^\times$, hence  $\gamma^{-n} k \gamma^n\in \gamma^{j-2n} N\cdot K'$ and $f_K(\gamma^{-n} k \gamma^n)=(-p^3)^{j-2n}$.

Therefore $\displaystyle \frac{1}{\vol(K')} \int\limits_{K'} f_K(\gamma^{-n} k \gamma^n) dk= \frac{1}{[K':K'_{2n}]} + p^{-6n} \sum_{j=0}^{2n-1} (-1)^j p^{3j}
\frac{1}{[K':K^{\prime\times}_j]}=$

\[  =c_0 \cdot p^{-6n} \left(  p^{2n} +c_0^{-1} -  \sum_{i=1}^{n}  p^{6i-3}(p^{-4i+3}- p^{-4i})  + \sum_{i=0}^{n-1} p^{6i}(p^{-4i} - p^{-4i-1}) \right)=
 p^{-6n} (1+ c_0).\qedhere \]
\end{proof}

\begin{rem} Alternatively,  one could use MacDonald's formula for zonal spherical functions  to see that $\pi_2$ does not admit  non-zero  vectors fixed by  the  hyperspecial maximal   open compact subgroup $K^\circ$. Indeed,  \cite[\S 5.5]{haines-kottwitz-prasad}
allows one to write $\Phi^{K^\circ}_n$ from \eqref{Phi-n},  up to a non-zero constant, as
\[   \Gamma_{\mu^{}_0} \cdot \mu^{}_0(\gamma^{-n}) +\Gamma_{\mu_0^{\sw}}\cdot \mu_0^{\sw}(\gamma^{-n}), \text{ with } \Gamma_\nu  =  \frac{1-p^2\cdot\nu(\gamma^{-2})}{1-\nu(\gamma^{-2})}\cdot
\frac{1-p^2\cdot\nu(\gamma^{-1})}{1-\nu(\gamma^{-1})}, \]
where the two factors in $\Gamma_\nu$  correspond respectively to the positive roots $\zeta$ and  $2\zeta$ of $G$ (see \S\ref{tree}).  As $\mu^{}_0(\gamma^{-1})=\mu_0^{\sw}(\gamma)=-p^{-1}$, one has  $\Gamma_{\mu^{}_0}=0\ne \Gamma_{\mu_0^{\sw}}$, hence the sequence $(\Phi^{K^\circ}_n)_{n\geqslant 0}$ is not $\rL^2$. This also shows, in passing, that $\pi_n$ is not square integrable.
\end{rem}

As a consequence  we obtain the following lower bound, in the unramified case.

\begin{cor}  \label{cor:deep-iwahori-inv}  For  $r\geqslant 1$ and for $\pi_n\in\Pi(\lambda_0)=\Pi(\lambda_0, \mathbf{1})$, one has $\dim\left(\pi_n^{I_{2r}}\right)\geqslant r+1$.
\end{cor}
\begin{proof} 
Let $f$ be a basis of $\pi_n^{K^\circ}$ (see Proposition~\ref{prop:inert-inv}). 
For all $r\in\Z$, $\pi_n$ contains a (unique) line fixed by  $\gamma^r K^\circ\gamma^{-r}$ (having  
$\gamma^r\cdot f_{K^\circ}$ as basis) and moreover, the stabilizer in $G$ of that line is $\gamma^r K^\circ\gamma^{-r}$.
We claim that the vectors  $f, \gamma\cdot f, \dots, \gamma^r\cdot f\in \pi_n$ are linearly independent. Indeed, if $f$ was a linear combination of the remaining $r$ vectors then it  would be fixed by $\displaystyle\bigcap_{1\leqslant j \leqslant r} \gamma^i K^\circ \gamma^{-i}$  which is {\it not} contained in $K^\circ$.  As any of  these $(r+1)$ vectors is  fixed by $I_{2r}$,  the claim follows. 

A similar claim  holds  for  $\pi_2$ and can be proven exactly as above, using $K'$ instead of $K^\circ$. 
 \end{proof}

While  the unramified $A$-packet $\Pi(\lambda_0)$ would  suffice for our global applications when $D\equiv 3\pmod{8}$, the case of  discriminants $D\equiv 7\pmod{8}$ would  require the use of  certain tamely ramified $A$-packets $\Pi(\lambda)$ and providing 
explicit levels for them is the object of the next subsection.

\subsection{Intertwinings and  exponential sums} \label{intertwining}

Our arithmetic applications will require to show  existence of non-zero $I_r$-invariants, for some $r\in \Z_{\geqslant 1}$,  in certain ramified $A$-packets. This is delicate because of the lack of a new-vector theory for non-generic representations $\pi$ (see Remark~\ref{new-vector}). Also Casselman's result  asserting that $\pi^{I_r}$ surjects onto $\pi_N^{T\cap I_r}$ is inconclusive here  as the latter vanishes, 
in contrast with the situation in \cite{dimitrov-ramakrishnan} where the open compact  is a pro-$p$-Iwahori subgroup of a sufficiently deep level. We will instead employ  explicit methods and prove  in Theorem~\ref{inert-inv-ram} that $\pi_n$  has non-zero invariants  by $K_T$, which contains a conjugate of $I_3$. It is relatively straightforward to determine all such  vectors  in the full induced representation, but it becomes a thorny comptutation to find a non-square integrable matrix coefficient as in \S\ref{iwahori-inv}. By another result of Casselman, matrix coefficients can be expressed in terms of the corresponding ones in the Jacquet module, given here by an explicit  scalar product on $\pi_N\simeq\C^2$.  Making this actually work requires non-vanishing of the second coordinate of the image of a $K_T$-invariant  vector under the Jacquet functor which, once  verified,  leads directly to the result we seek. 

We  reduce the  computation of the Jacquet functor  to a precise statement about the intertwining operator at the level of finite reductive groups. The latter involves showing non-vanishing of some explicit exponential sums, bringing out the arithmetic nature of the problem. Although these sums seem extremely hard 
to  compute  individually, we manage to deduce it by precisely computing a suitable average of such sums corresponding to the trace of  the finite intertwining operator.

\begin{thm}\label{inert-inv-ram} Assume  that $p$  odd  and  $E/\Q_p$  unramified. Fix a character $\chi:E^1\twoheadrightarrow \F^1_{p^2}\to \C^\times$ such that $\chi^3\ne \mathbf{1}$ and let  $\lambda=\lambda_0\cdot \chi^{}_E$. Letting $\pi_n$ denote the non-tempered representation of the Arthur packet $\Pi(\lambda)$, one has $\dim\pi_n^{K_T}=1$.
\end{thm}

\begin{lemma} \label{KT-inv} Consider the character $\mu_\lambda$ from \eqref{eq:mu} attached to $\lambda$ and trivial on $E^1$. 
One has
\[B\backslash G /K_T= \overline{B}\backslash \overline{G} /\overline{T}= \left\{\mathbf{1}_3, \sw,
\sigma_\infty=\left(\begin{smallmatrix} 1 & 0 & 0\\ 0  & 1 &0 \\ \xi   & 0 & 1 \end{smallmatrix}\right),
\sigma_y=\left(\begin{smallmatrix} 1 & 0 & 0\\ -1  & 1 &0 \\ \xi y-\frac{1}{2}  & 1 & 1 \end{smallmatrix}\right)\Big{|} y\in \F_p\right\}. \]

For any non-trivial $\chi$,  the $K_T$-invariant functions in $\Ind_B^G(\mu_\lambda)$ are supported by the double cosets of $\{\sigma_y \vert y\in \F_p\}$ and,  in addition when $\chi$ is quadratic, by the double coset of $\sigma_\infty$.
\end{lemma}

\begin{proof} By the Iwasawa decomposition, restriction to $K$ yields an isomorphism  
$\Ind_B^G(\mu_\lambda)\xrightarrow{\sim} \Ind_{B\cap K}^K(\mu_\lambda)$. As $\mu_{\lambda|B\cap K}$ factors through a character 
$\overline{\mu}_\lambda$ of $\overline{B}$, the space of $K_1$-invariants in  $\Ind_{B\cap K}^K(\mu_\lambda)$ is naturally  identified  with  $\Ind_{\overline{B}}^{\overline{G}}(\overline{\mu}_\lambda)$, and therefore its  subspace of 
$K_T$-invariants  is given by $\Ind_{\overline{B}}^{\overline{G}}(\overline{\mu}_\lambda)^{\overline{T}}$. 
By definition, a double coset $\overline{B}\sigma\overline{T}$ supports a non-zero function contributing to the latter space if, and only if, 
$\overline{\mu}_\lambda$ is trivial on $\overline{B}\cap \sigma \overline{T}\sigma^{-1}$
(this is in accordance with Mackey's Theorem for finite groups). 
As $\overline{\mu}_\lambda$ is ramified this is never the case for $\sigma=\mathbf{1}_3$ nor for $\sigma=\sw$, while 
a matrix computation shows that $\sigma=\sigma_\infty$ works if, and only if, $\chi$ is quadratic.
For $y\in \F_p$ and  $k\in K_T$, writting $(k\bmod{p})= \diag \left( \bar\alpha, \beta,\alpha^{-1} \right)$, we find 
 \[(\sigma_y k \sigma_y^{-1}\bmod{p})=\left(\begin{smallmatrix} \bar\alpha & 0 & 0\\ \beta- \bar\alpha & \beta & 0 \\  * & \beta - \alpha^{-1} & \alpha^{-1} \end{smallmatrix}\right)\] 
 ensuring the triviality of $\mu_\lambda$ on $B\cap \sigma_y K_T \sigma_y^{-1}$, as in forces $\bar\alpha\equiv \beta \pmod{p}$ on such elements. 
 \end{proof}

 Recall that $\pi_n$ is the Langlands quotient of $\Ind_B^G(\mu_\lambda)$  whose other  Jordan-H\"older constituent is the discrete series $\pi_2$. 
As taking invariants by an open compact subgroup is an exact functor in the category of admissible representations,  Lemma~\ref{KT-inv} implies that    
\begin{equation}\label{eq:dim-sum}
\dim \pi_n^{K_T}+ \dim \pi_2^{K_T}=\begin{cases} p & \text{ if } \chi^2\ne \mathbf{1},\\p+1 & \text{if } \chi \text{ quadratic}.\end{cases}   
\end{equation}
\begin{rem} \label{new-vector}
To show  $\pi_n^{K_T}\ne \{0\}$ one could  instead try  computing $\dim \pi_2^{K_T}$.  
Miyauchi's theory \cite{miyauchi-U21, miyauchi-L-factor} of conductors with respect to the  paramodular groups $K'_r=\left(\begin{smallmatrix}\cO^\times & \cO &p^{-r}\cO\\p^r\cO &\cO^\times &\cO\\ p^r\cO & p^r\cO &\cO^\times\end{smallmatrix} \right)\cap G$ asserts that the level of the generic $\pi_2$ is given by its conductor, while 
$\pi_c$ and the ramified $\pi_n$  have no level, {\it i.e.} they have no non-zero invariants by $K'_r$ for any $r$.   When $\chi$ is the quadratic character,  the  $L$-parameter $(\lambda \otimes \mathrm{St})  \oplus \mathbf{1}$ has conductor $2$, therefore the generic member  $\pi_2$ in this $L$-packet has $1$-dimensional invariants by $K'_2$, hence also by  $\gamma^{-1}K'_2\gamma= \left(\begin{smallmatrix}\cO^\times & \cP &\cO\\ \cP &\cO^\times &\cP\\ \cO & \cP &\cO^\times\end{smallmatrix} \right)\cap G\supset K_T$. For other ramified $\chi$'s,  $\pi_2$ has  an invariant line by $K'_3$, having the   same volume as $K_T$ but  not conjugate to it, thus not settling the non-vanishing of $\pi_2^{K_T}$, let alone computing its dimension. 
\end{rem}

In order to show that $\pi_n$ admits non-zero $K_T$-invariants, it would be enough to consider the image of a $K_T$-invariant vector from Lemma~\ref{KT-inv} under the  Jacquet functor and show that its second coordinate does not vanish. As well known, this functor  is  given by 
\[\Ind_B^G(\mu_\lambda) \longrightarrow  \C\cdot \mu_\lambda \delta^{1/2} \oplus  \C\cdot\mu_\lambda^{\sw}\delta^{1/2}, \quad  f\mapsto (f(1),(Mf)(1))\]
where the standard intertwining operator  $M:\Ind_B^G(\mu_\lambda) \to \Ind_B^G(\mu_\lambda^{\sw})$ is defined via analytic continuation, as follows. 
For $s\in \C$, letting $\mu_{\lambda,s}=\mu_\lambda\cdot \delta^{s/2}$, the  intertwining operator 
\[M_s:\Ind_B^G(\mu_{\lambda,s}) \to \Ind_B^G(\mu_{\lambda,s}^{\sw}), \quad f_s\mapsto \int\limits_N f_s(\sw n\cdot)dn\]
is absolutely convergent for  $\Re(s)\gg 0$ and  $G$-equivariant. 
 Moreover, for  any section $f_s\in \Ind_B^G(\mu_{\lambda,s})$ such that for all $g\in G$ the function $f_s(g)$ is analytic in $s\in \C$, the function $(M_s(f_s))(g)$, a priori only defined for $\Re(s)\gg 0$, is a rational function in $p^{-s}$, hence extends to a meromorphic function on all of $\C$ with only possibly a finite number of poles independent of $f$ and of $g$. In fact, it continues as an intertwining operator, {\it i.e.} 
\[  (M_s(f_s))(gg')= (M_s(f_s(\cdot g'))(g), \text{ for all } g,g'\in G. \]
We refer to \cite[\S 1]{arthur-intertwining} for more detail and proofs.   We will  computing explicitly the right hand side~of
\[M(f_0)=\lim_{s\to 0}  M_s(f_{s})\]
 as a rational function in $p^{-s}$, simultaneously justifying the absence of pole at $s=0$. 

 For any $y\in \F_p$, Lemma \ref{KT-inv} yields a unique function   $f_y\in\left(\Ind_B^G(\mu_\lambda)\right)^{K_T}$ supported on  $B\sigma_y K_T$ and normalized by $f_y(\sigma_y)=1$. In addition, when  $\chi$ is quadratic, we let  $f_\infty$ be the unique function supported on  
  $B\sigma_\infty K_T$ and normalized by $f_\infty(\sigma_\infty)=1$.
 Consider a $K^\circ$-flat section  $f_{y,s}$ passing through $f_{y,0}=f_{y}$. 
 
\begin{prop} \label{Ms-fy} 
For all $y,y'\in \F_p\cup\{\infty\}$ and $\Re(s)\gg 0$, we have
\[M_s(f_{y',s})(\sigma_y)=\chi(-1) \frac{ (1-p)p^{-3-2s}}{1-p^{-1-2s}}\delta_{y,y'}+p^{-3} \sum_{[z,x]\in N(\F_p) }
f_{y',s}\left(\sw [z,x]\sigma_y\right).\]
\end{prop}

\begin{proof}
Recall that $[z,x]= \left(\begin{smallmatrix}1  & -\bar{z}& u  \\ &1 & z\\  & & 1\end{smallmatrix} \right)$, where
$u=\xi\cdot x - \tfrac{z\bar z}{2} $, so  that  $u+\bar u + z \bar z =0$. 
 We first   compute the intertwining  integral over 
\[N(\Q_p){\smallsetminus} N(\Z_p)=\bigcup_{m\geqslant 1 \text{ odd}}[p^{\frac{1-m}{2}}\cO , p^{-m} \Z_p^\times]
\bigcup_{m\geqslant 2 \text{ even}}[p^{\frac{2-m}{2}}\cO , p^{-m} \Z_p^\times]
\bigcup_{m\geqslant 2 \text{ even}}[p^{-\frac{m}{2}}\cO^\times, p^{-m} \Z_p],\]
a key observation being that $[z,x]\in N(\Z_p)$ if and only if  $u \in \cO$. It follows that 
for  $[z,x]\in N(\Q_p){\smallsetminus} N(\Z_p)$  both $1/u$ and  $-\bar{z}/u$ belong to $\cP$, yielding  an explicit Iwasawa decomposition:   
\[\sw n= \sw [z,x]=\begin{pmatrix} &  & 1 \\  & 1 & z \\ 1  & -\bar{z} & u \end{pmatrix}=
 \begin{pmatrix}\bar u^{-1}  & \bar{z}/u & 1 \\  & -\bar u/u & z \\  & & u \end{pmatrix}
\begin{pmatrix} 1&  &  \\   z/\bar{u} & 1 & \\ 1/u  & -\bar{z}/u & 1 \end{pmatrix}.\]
As the last term belongs to the principal congruence subgroup $K_1$ which is normalized by $\sigma_y\in K^\circ$ and contained in $K_T$, 
we deduce that $f_{y',s}(\sw n\sigma_y)=\begin{cases}\mu_{\lambda,s}(\bar u^{-1}, -\bar{u}/u, u)& \text{ if } y'=y, \\ 0  & \text{ if } y'\ne y.\end{cases}$ 
Hence for $u\in p^{-m}\cO^\times$, we have 
$f_{y',s}(\sw n\sigma_y)=\mu_{\lambda,s}(\gamma)^{-m}\chi^3_E(\bar u^{-1})=(-1)^m \cdot p^{-(3+2s)m}\chi^3_E(u)$.

For $m$ odd, noticing further that   $\chi^{}_E(u)=\chi^{}_E(\xi)=\chi(-1)$, we find
 \[\int\limits_{[p^{(1-m)/2}\cO , p^{-m} \Z_p^\times]}  f_{y,s}(\sw n\sigma_y)dn
 =-\chi(-1)\cdot p^{-(3+2s)m}    p^{m-1} (p^m-p^{m-1}), \text{ hence} 
 \]
 \begin{align}\label{eq:int1}
 \sum\limits_{m\geqslant 1 \text{ odd}}  \int\limits_{[p^{(1-m)/2}\cO , p^{-m} \Z_p^\times]}  f_{y,s}(\sw n\sigma_y)dn 
 =  \frac{-\chi(-1) (p^{-2-2s}-p^{-3-2s}) }{1-p^{-2-4s}}. 
 \end{align}

 We similarly  find 
  \begin{align}\label{eq:int2}
 \sum\limits_{m\geqslant 2 \text{ even}} \int\limits_{[p^{(2-m)/2}\cO , p^{-m} \Z_p^\times]}  f_{y,s}(\sw n\sigma_y)dn
 =  \frac{\chi(-1) (p^{-4-4s}-p^{-5-4s}) }{1-p^{-2-4s}}.
  \end{align}
 Finally, for $m\geqslant 2$ even we compute 
 \[\int\limits_{[p^{-\frac{m}{2}}\cO^\times, p^{-m} \Z_p]}  f_{y,s}(\sw n\sigma_y)dn
 =p^{-(3+2s)m}   \cdot p^{m-2} p^{m-1} \sum\limits_{[z_0,x_0]\in\F_{p^2}^\times \times  \F_p} \chi^{3}_E(z_0\bar{z}_0-2\xi x_0).  
 \]
 As $\chi^3\ne \mathbf{1}$, 
 $\sum\limits_{(z_0,x_0)\in\F_{p^2}^\times \times  \F_p} \chi^3\left(\frac{z_0\bar{z}_0-2\xi x_0}{z_0\bar{z}_0+2\xi x_0}\right)=(p^2-1)
 \sum\limits_{x_0 \in \F_p} \chi^3\left(\frac{1-\xi x_0}{1+\xi x_0}\right)= -\chi(-1)(p^2-1)$ yields
  \begin{align}\label{eq:int3}
\sum\limits_{m\geqslant 2 \text{ even}} \int\limits_{[p^{-\frac{m}{2}}\cO^\times, p^{-m} \Z_p]}  f_{y',s}(\sw n\sigma_y)dn= 
   \frac{-\chi(-1) (p^{-3-4s}-p^{-5-4s}) }{1-p^{-2-4s}}. 
  \end{align}
Putting \eqref{eq:int1}, \eqref{eq:int2} and \eqref{eq:int3} together, we find 
$\int\limits_{N(\Q_p){\smallsetminus} N(\Z_p)}f_{y,s}(\sw n\sigma_y)dn= -\chi(-1) 
\frac{(p^{-2-2s}-p^{-3-2s}) }{1-p^{-1-2s}}$. 

Finally, as $f_{y',s}$ is invariant by $K_T\supset K_1\supset N(p\Z_p)$, we have 
\[\int\limits_{N(\Z_p)}f_{y',s}(\sw n\sigma_y)dn=p^{-3} \sum_{[z,x]\in N(\Z_p)/N(p\Z_p) }
f_{y',s}\left(\sw [z,x]\sigma_y\right).\qedhere \]
\end{proof} 

\noindent The coefficients of the intertwining automorphism of $\Ind_{\overline{B}}^{\overline{G}}(\overline{\mu}_\lambda)^{\overline{T}}$ in the basis  $(f_{y})_y$ are given by
\[M_{y,y'}=\sum_{[z,x]\in \F_{p^2}\times\F_p }f_{y',0}\left(\sw [z,x]\sigma_y\right).\]

\begin{cor} \label{cor:M0-eigenvalues} Letting $\mathscr{M}_0$ (resp. $\mathscr{M}_0^{\sw}$) denote  the matrices  of $M_0$ (resp. $M_0^{\sw}$) in the bases $(f_{y})_y$ and $(f^{\sw}_{y})_y$, we have 
\[\mathscr{M}:=(M_{y,y'})_{(y,y')}=\chi(-1) p \cdot I_p + p^3\cdot\mathscr{M}_0= -\chi(-1) p^2  \cdot I_p + p^3\cdot\mathscr{M}_0^{\sw}.\]
\end{cor}
\begin{proof} The first equality is obtained by letting $s$ go to $0$ in Proposition~\ref{Ms-fy}, while the second one is obtained by applying the same process to
 $\mu_\lambda^{\sw}$, noting that   $\mu_{\lambda,s}^{\sw}=\mu_{\lambda,s-1}$, in particular $\overline{\mu}_\lambda=\overline{\mu}_\lambda^{\sw}$. 
 \end{proof}

\begin{lemma} \label{lem:M0-trace}
For all $y\in \F_p$, one has $\sum\limits_{y'\in \F_p}M_{y,y'}=\sum\limits_{y\in \F_p}M_{y,y'}=\chi(-1) p$. More importantly 
\begin{align}\label{eq:trace-M}
\Tr(\mathscr{M})=\begin{cases} -\chi(-1)p & \text{ if } \chi^2\ne \mathbf{1} \quad(\text{and} \chi^3\ne \mathbf{1}), \\0 &  \text{ if }  \chi \text{ quadratic}. \end{cases} 
\end{align}
\end{lemma}

\begin{proof} As $f_{y'}(\sw\sigma_y)= 0$, we can assume that $u\ne 0$. 
Letting $v=1/u$ and $w=1-\bar z/u$,  we have   
\[\mathsf{w}[z,x]\sigma_y=\begin{pmatrix}-\tfrac12+\xi y & 1 & 1 \\ -1-\tfrac{z}{2}+\xi yz & 1+z & z \\ 1+\bar z -\tfrac{u}{2} +\xi yu & u- \bar z & u \end{pmatrix}=
\begin{pmatrix}\bar u^{-1}  & * & * \\  & -\bar u/u & * \\  & & u \end{pmatrix}
\begin{pmatrix}1  &  &  \\ -\bar w & 1&  \\  v-w +\tfrac12 +\xi y& w & 1 \end{pmatrix}
 \]

Hence  $f_{y'}(\sw [z,x]\sigma_y)\ne 0$ if and only if $v-w +\tfrac12 +\xi y=w\bar w (-\tfrac12+\xi y')\ne 0$,  {\it i.e.},
\begin{equation} \label{eq:yy-prime} 
v- \bar v = 2 \xi (w \bar w y'- y)+w- \bar w  \, \text{  and } w\ne 0
\end{equation}
in which case it equals $\chi^{}_E(u^3w )=\chi^{}_E( w)\cdot  \chi^3_E(\bar v)$. Noting $(w-1)(\bar w-1) + v+\bar v=0$, so that for a fixed $w\ne 1$, the ratio 
$\tfrac{\bar v}{v}=\tfrac{u}{\bar u}=\frac{\xi x -(w-1)(\bar w-1)/2}{-\xi x -(w-1)(\bar w-1)/2}$ describes 
$\F^1_{p^2}{\smallsetminus} \{-1\}$, we find
\begin{align*} 
&\sum_{y\in \F_p} M_{y,y'}=\sum_{y'\in \F_p} M_{y,y'}=\\
&=\sum_{w\ne 1} \chi^{}_E( w)  \sum_{x\in \F_p}  \chi^{3}_E\left(\xi x -(w-1)(\bar w-1)/2\right)+\sum_{x\in \F_p^\times} \chi^3(-1)\\
&=\sum_{w\ne 1} \chi^{}_E(w)(-\chi(-1))+ (p-1) \chi(-1) = p\chi(-1).
\end{align*} 

To compute $\Tr(\mathscr{M})$ we observe that when $w\notin \F^1_{p^2}$,   \eqref{eq:yy-prime} imposes no condition on $v$, as we are summing over all 
$y\in \F_p$, and that  when $w\in \F^1_{p^2}$, it  imposes no condition on $y$ while determines uniquely 
$v=\tfrac{w- \bar w }{2}-\tfrac{(w-1)(\bar w-1)}{2}=w-1$ . As $v\ne 0$ implies $w\ne 1$, we find 
\begin{align*} 
&\sum_{y\in \F_p} M_{y,y} =- \chi(-1) \sum_{w\in \F_{p^2}^\times{\smallsetminus} \F^1_{p^2}} \chi^{}_E(w)  
+p \sum_{w\in  \F^1_{p^2}{\smallsetminus} \{1\}} \chi^{}_E(w (\bar w-1)^3) \\
&=\chi(-1) \sum_{w\in \F^1_{p^2}} \chi^2(w) + p \sum_{w\in  \F^1_{p^2}{\smallsetminus} \{1\}} \chi(-\bar w)=
\begin{cases} -p \chi(-1) & \text{ if } \chi^2\ne \mathbf{1}, \\ \chi(-1) &  \text{ if }  \chi \text{ quadratic}. \end{cases} 
\end{align*}  

Finally, for $\chi$ quadratic,  letting $w=-\tfrac{\bar z}{u}$, we find 
$\mathsf{w}[z,x]\sigma_\infty=
\left(\begin{smallmatrix}\bar u^{-1}  & * & * \\  & -\bar u/u & * \\  & & u \end{smallmatrix}\right)
\left(\begin{smallmatrix}1  &  &  \\ -\bar w & 1&  \\  u^{-1}  +\xi & w & 1 \end{smallmatrix}\right)$. 
This implies that $f_{\infty}(\sw [z,x]\sigma_\infty)=0$ unless $w=z=0$ in which case it equals $\chi^{}_E(\xi y)$, 
with $y\bar y= 1+\xi^{-2} x^{-1}\ne 0$.  Then one finds $M_{\infty,\infty}=\sum\limits_{x\in \F_p^\times }f_{\infty}\left(\sw [0,x]\sigma_\infty\right)=-\chi(-1)$.
Combining with the second case of the above  equation,  one finds $\Tr(\mathscr{M})=M_{\infty,\infty}+ \sum\limits_{y\in \F_p} M_{y,y}=0$. 
\end{proof}

\begin{proof}[Proof of Theorem~\ref{inert-inv-ram}] 
As  $\overline{\mu}_\lambda=\overline{\mu}_\lambda^{\sw}$,  Frobenius reciprocity yields 
 $\End_{\overline{G}}\left(\Ind_{\overline{B}}^{\overline{G}}(\overline{\mu}_\lambda)\right)=\C^2$. Hence 
 $\Ind_{\overline{B}}^{\overline{G}}(\overline{\mu}_\lambda)$ is a direct sum of two irreducible distinct $\overline{G}$-representations, denoted  $\overline{\pi}_n$ and $\overline{\pi}_2$, on each of which the intertwining endomorphism  $\mathscr{M}$ acts scalarly by Schur's Lemma. 
 
 Furthermore, as $\overline{\pi}_n^{\overline{T}}$ and $\overline{\pi}_2^{\overline{T}}$  are naturally  identified with $\pi_n^{K_T}=\ker(\mathscr{M}_0^{\sw})$ and  
 $\pi_2^{K_T}=\ker(\mathscr{M}_0)$ (see \eqref{induction}),  it follows from  Corollary~\ref{cor:M0-eigenvalues}  
 that the above mentioned eigenvalues are respectively $-\chi(-1) p^2$  and $\chi(-1) p \cdot I_p$. 
Formula \eqref{eq:trace-M} can thus be rewritten as 
\[ p\cdot \dim\pi_n^{K_T} - \dim\pi_2^{K_T} =-\frac{\chi(-1)}{p}\Tr(\mathscr{M})=
\begin{cases} 1 & \text{ if } \chi^2\ne \mathbf{1} \quad(\text{and} \chi^3\ne \mathbf{1}), \\0 &  \text{ if }  \chi \text{ quadratic}, \end{cases} \]
which combined with \eqref{eq:dim-sum} yields $\dim\pi_n^{K_T}=1$. 
\end{proof}

\section{Galois representations for  $3$-folds of Picard type}\label{MT-conj}

From this point onwards, we will use  global notations from the introduction.  The  local results of the previous section can be applied to the completion  $E$ of $M$ at any prime number which does not split in that field.
We denote by $\A_f$  the ring of finite adeles of $\Q$, so that $\A=\A_f\times \R$.

\subsection{Abelian $3$-folds and Shimura surfaces of Picard type} \label{sec:moduli}
Let $k$ be any field containing $M$.
 Consider  an abelian $3$-fold  $A/k$   together with an injection
 $\iota^0: M \hookrightarrow \mathrm{End}^0(A/k)=\mathrm{End}(A/k)\otimes\Q$, or equivalently with an injection $\iota$
 of an  order of $M$  into $\mathrm{End}(A/k)$, the most important for us case being  when  $\iota^0$  comes from $\iota: \cO_M\hookrightarrow
\mathrm{End}(A/k)$.

A polarization on $A$ is an isogeny $\theta: A\to A^\vee$, where $A^\vee$ denotes the dual abelian variety.
By positivity, since $k$ is a field, the Rosati involution induced by $\theta$ on $\iota^0(M)$ is given by  the complex conjugation (see \cite[\S 21]{mumford}). A polarization is called  principal, if it is an isomorphism, and it can always be acquired over a finite extension of  $k$. 

By an abelian $3$-fold of {\it Picard type} over $k$ we will mean a principally polarized abelian 
variety over $k$ of dimension $3$   having multiplication by $\cO_M$ defined over $k$. 

The action of $M$ splits the $3$-dimensional $k$-vector space $\Lie(A/k)$  in a direct sum of two sub-spaces:  one on which the actions of $M$ and $k$ agree, and  one on which they differ by the complex conjugation, the pair of their dimension being called the signature. 

To define a level structure on $A$ we need to  consider its Tate module.
Given a finite place $v$ of $M$,  the $v$-adic Tate module $T_v A=\varprojlim\limits_r A[v^r]$ of $A$  is free of rank $3$  over $\cO_v$.
Denote $V_v A=M_v \otimes_{\cO_v}T_v A $.
One also considers the adelic Tate module
 \[V_f A=\Q \otimes_{\Z}  \varprojlim\limits_{n}A[n],\]
 which is free of rank $3$ over  $\A_{M,f}$.
Given a polarization $\theta: A\to A^\vee$, the Weil pairing endows  $V_f A$ with a non-degenerate skew-hermitian form,
{\it i.e.}, a non-degenerate alternating pairing
\[\langle \cdot, \cdot \rangle_A:V_f A\times V_f A\to \A_f\]
such that $\langle a\cdot v, v'  \rangle_A=\langle  v, \bar a\cdot v'  \rangle_A$ for all $a\in M$.
If $\theta$ is principlal, then   $\langle \cdot, \cdot \rangle_A$ is a perfect pairing.

Let $(V,\langle \cdot, \cdot \rangle)$ be a $3$-dimensional (non-degenerate) hermitian space over $M$.
The unitary similitude  group $\wG=\mathrm{GU}(V)$ is the reductive group over $\Q$ characterized by the property that
\[\wG(R)= \{ g\in \GL(V\otimes_\Q R) \mid \forall v,v'\in V\otimes_\Q  R,  \langle g(v), g(v') \rangle=\nu(g) \langle v, v' \rangle\},\]
for any $\Q$-algebra $R$,  where $\nu: \mathrm{GU}(V)\to \G_{m,\Q}$ is a homomorphism whose kernel is the  unitary group $ G= \rU(V)$. Since  any hermitian form in $3$ variables over a non-archimedean local field is isotropic, the group $\wG (\Q_p)$ is unique up to isomorphism, while at infinity $\langle \cdot, \cdot \rangle$ is uniquely determined by its signature,  hence there are  two possibilities for $\wG(\R)$ (as  opposite signatures define isomorphic groups).  Hasse's Principle applied to the semi-simple simply connected derived group $G^1=\mathrm{SU}(V)$ implies then   that, up to an isomorphism, there exists a unique quasi-split unitary  group, denoted $\wG=\mathrm{GU}(2,1)$, and a unique definite  unitary  group, denoted $\mathrm{GU}(3)$.

We will now  introduce the Shimura variety for  $\wG$
 represented by the matrix $\left(\begin{smallmatrix} & & 1 \\ & \sqrt{-D} & \\ -1 & & \end{smallmatrix} \right)$.
The   homomorphism  
\[\tilde h: \mathrm{Res}^\C_\R \G_{m,\R} \to \wG_\R \text{ , } z\mapsto \begin{pmatrix} \Re(z) & 0 &  \Im(z) \\
0 & z & 0 \\ -\Im(z)  & 0 &  \Re(z) \end{pmatrix}\]
of  $\R$-algebraic   groups  satisfies Deligne's axioms \cite{deligne-ShVar} for a  Shimura variety, hence for any open compact subgroups  $\wK$  of  $\wG(\A_f)$ one can consider the Shimura surface
\[ Y_{\wK}(\C) =\wG(\Q)\backslash \left(\mathcal{H} \times \wG(\A_f)/\wK\right), \]
where  $\mathcal{H}$ is identified  with the $\wG(\R)$-conjugacy classes of $\tilde h$.
By a  fundamental result of  Shimura, $Y_{\wK}$ admits  a canonical model over the reflex field $M$. The connected components of  $Y_{\wK}$ are Picard modular surfaces.

For the anisotropic form $\mathrm{GU}(3)$, one can analogously define Shimura  sets which are finite and therefore will not alter the uniformity of our results in \S\ref{s:uniform}.

The Shimura surfaces of Picard type are coarse moduli spaces of abelian $3$-folds. Namely,
$Y_{\wK}(\C)$ is in bijection with isogeny classes of $(A,\iota^0,\theta, \eta\cdot \wK )$, where
$(A,\iota^0, \theta)$ is an  abelian $3$-fold over $\C$ as above, and $\eta: \A_f\otimes_{\Q}V \xrightarrow{\sim} V_f A$ is an isomorphism sending $\langle \cdot, \cdot \rangle_A$ to a 
$\A_f^\times$-multiple of  $\langle \cdot, \cdot \rangle_V$.
Note that the usual $\Q^\times$-multiple condition is automatically satisfied as we are in the type C case.
When $\wK^\circ$  is the standard maximal open compact subgroup of $\wG(\A_f)$,  the points of  $Y_{\wK^\circ}(k)$ correspond to  isomorphism classes of abelian $3$-folds of Picard type.

Let  $(A,\iota, \theta)$ be an  abelian $3$-fold of Picard type over a number field $k$. 
 The action of $\Gal_k$ on the adelic Tate module $V_f A$, together with the choice of $\eta$ as above,  yields  a continuous homomorphism:
\begin{align}\label{eq:galois-action}
 \rho_{A,f}: \Gal_k\longrightarrow \wG(\A_f). 
 \end{align} 
Moreover, the point $(A,\iota,\theta, \eta \cdot \wK )$ on
$Y_{\wK}(\C)$ is defined over $k$ if, and only if, $\rho_{A,f}(\Gal_k)\subset \wK$.

As for each $\ell$ the principal polarization $\theta$  induces a perfect pairing on $T_\ell A\simeq \rH_1(A(\C), \Z_\ell)$, it follows that  one can choose $\eta$ so that  $\rho_{A,f}(\Gal_k)\subset \wK^\circ$, {\it i.e.}, $(A,\iota, \theta, \eta\cdot \wK^\circ)$ defines a $k$-rational point on $Y_{\wK^\circ}$.
 By the  Brauer--Nesbitt Theorem the  semi-simplification $\overline\rho_{A,\ell}$ of the composition of $\rho_{A,\ell}: \Gal_k\longrightarrow \wK_\ell^\circ$ with the natural  surjection of  $\wK_\ell^\circ$ onto its reductive quotient $\overline{G}_\ell^\circ$
 is uniquely determined by its characteristic polynomial. Moreover,  $\overline\rho_{A,\ell}$ is absolutely irreducible if, and only if, the  self-dual lattice   fixed by $\wK_\ell^\circ$ is the  unique, up to homothecy, $\rho_\ell(\Gal_k)$-stable lattice.

  In analogy with   the  index $2$ Gross subgroups $\wK''_p$  of maximal compacts  $K^\circ_p$ at ramified primes introduced in \S\ref{K-second}, we consider the open compact subgroup
\begin{equation}\label{K-double-prime}
\wK''= \wK''_D \prod_{p\nmid D} \wK_p^\circ \subset \wG(\A_f),
\end{equation}
where $\wK''_D$ is  defined as the kernel of the composed homomorphism
\begin{equation}\label{KD-double-prime}
 \prod_{p\mid D}  \wK_p^\circ  \twoheadrightarrow \prod_{p\mid D}  \wK_p^\circ/\wK''_p= \prod_{p\mid D} \{\pm 1\}
 \xrightarrow{\Pi} \{\pm 1\}.
 \end{equation}

Let $(A,\iota,\theta)$ be a principally polarized  abelian $3$-fold of Picard type over $k$.
For $v$ the prime of $M$ above  $p\mid D$, the action of  $\Gal_k$  on $A[v]$
 factors through $\mathrm{GO}(3,\F_p)$.
 Using the exceptional isomorphism
 $\mathrm{PGO}(3,\F_p)\xrightarrow{\sim} \mathrm{SO}(3,\F_p)\xrightarrow{\sim}  \PGL(2,\F_p)$, one  defines its projectivization
\begin{equation}\label{eq:star}
\widetilde\rho_{A,p}: \Gal_k \to \PGL(2,\F_p).
\end{equation}
Taking quotient by the unique index two subgroup $\PSL(2,\F_p)$ of $\PGL(2,\F_p)$ yields a canonical homomorphism
$\varepsilon_{A,p}:  \Gal_k \to\{ \pm 1 \}$ and we let $\varepsilon_{A,D}=\prod_{p \mid  D} \varepsilon_{A,p}:  \Gal_k \to\{ \pm 1 \}$.

  Using the  observation made  after \eqref{eq:galois-action}, the points  in  $Y_{\wK''}(k)$ corresponds precisely to an abelian $3$-fold $A$ over $k$ of Picard type having trivial $\varepsilon_{A,D}$.

\subsection{Complex reflexions and families of abelian threefolds of Picard type}\label{sec:abelian-family}
It is important to observe  that even though for $\wK\subset \wK^\circ$ each point of $Y_{\wK}(\C)$ is associated to an abelian $3$-fold of Picard type,  there does not exist such a family over the entire $Y_{\wK}(\C)$ unless there are no points with extra automorphisms, in which case $Y_{\wK}(\C)$ would be a fine moduli space. In our cases of interest  $\wK$ is not neat, and therefore  $Y_{\wK}$ is not a fine moduli space.  
In this subsection we prove that there exists a natural family of  abelian $3$-folds of Picard type over  $Y_{\wK''}$ minus a finite number of points (see \eqref{K-double-prime}). This will be crucially used in the proof of Theorem~\ref{theoremB}.

Recall that  an element $\gamma$ of the discrete subgroup $\Gamma=\wG(\Q)\cap  \wK \wG(\R)$  has a fixed point in $\mathcal{H}$ if, and only if, $\gamma$ has finite order (this is because the stabilizers in $\wG(\R)$ of  points in $\mathcal{H}$ are maximal compact subgroups). 
Such a  $\gamma$ is called elliptic if it only fixes an isolated point in $\mathcal{H}$, otherwise it is called a complex reflexion. 
By \cite[Cor.~4.5.10]{holzapfel-book}, the set of singular points of $Y_\Gamma$ is finite and consists of isolated fixed points of elliptic 
elements of $\Gamma$, all defined over a finite extension of~$M$. Moreover, smooth points at which universal cover $\mathcal{H}\to Y_\Gamma$ is not \'etale are fixed points of complex reflexions, {\it i.e.}  order $2$ elements in $\Gamma$ fixing hyperbolic planes  (see  \cite[\S 4.5]{holzapfel-book}).

Given  a  geometrically connected component $Y''$ of $Y_{\wK''}\times_{M} k$ ,  we have $Y''(\C)=\Gamma''\backslash\mathcal{H}$ with $\Gamma''=\wG(\Q)\cap g_f\wK'' g_f^{-1}  \wG(\R)$ for some 
$g_f\in  \wG(\A_f)$.

\begin{lemma}\label{lem:reflexions}
The group $\Gamma''$ does not contain any complex reflexions.
\end{lemma}

\begin{proof} Consider a  complex reflexion $\gamma\in \Gamma''$ as an endomorphism of the Hermitian space $M^3$ having signature $(2,1)$. As the eigenspaces of $\gamma$ are mutually orthogonal, at most one such eigenspace could contain a negative line (corresponding to a point in $\mathcal{H}$).
This if $\gamma$ fixes more that one point of $\mathcal{H}$, it necessarily fixes a hyperbolic line in  $\mathcal{H}$. The corresponding endomorphism of  $M^3$ has a mutually orthogonal
 eigenplane and   eigenline  (both $M$-rational),  forcing the eigenvalues  to be in $\cO_M^\times=\{\pm 1\}$ (as $D>4$) and not all equal.   It follows, that for any $p\mid D$,  the image of $g_f^{-1}  \gamma  g_f\in\wK^\circ$  into the projectivization of the reductive quotient of $\wK_p^\circ$ belongs to the image under the adjoint isomorphism $\PGL(2,\F_p)  \xrightarrow{\sim}  \mathrm{PGO}(3,\F_p)$
 of an element represented by a matrix having both eigenvalues $1$ and $-1$. In particular its image in
 $\PGL(2,\F_p)/\PSL(2,\F_p)=\{\pm 1\}$ equals $\left(\tfrac{-1}{p}\right)$. As $\prod_{p\mid D} \left(\tfrac{-1}{p}\right)=  \left(\tfrac{-1}{D}\right)=-1$
 it follows from  \eqref{K-double-prime}  that $g_f^{-1}  \gamma  g_f\notin\wK''$, {\it i.e.}, $\gamma \notin \Gamma''$.
  \end{proof}
  \begin{prop}\label{ab-family}
There exists a finite extension $k$ of $M$ and a natural family of abelian $3$-folds of Picard type over $Y_{\wK''}\times_{M} k$ minus a finite number of $k$-rational elliptic points. 
\end{prop} 

\vspace{-10mm}

\begin{proof} 
We claim that there is a family of abelian $3$-folds of Picard type over any open subset $U$ of  $Y_{\wK''}$ which contains no point with a non-trivial stabilizer. 
  By   \cite[\S 2.3.4]{kwlan-book},  the moduli stack $S_{\wK''}$   associated to this problem is an algebraic stack (for the \'etale topology), locally of finite type over the base which we may take to $\Spec(M)$. Moreover, by 
  \cite[\S A.7.5]{kwlan-book}, there is a canonical surjective morphism $\phi$ from $S_{\wK''}$ to the associated coarse moduli space $[S_{\wK''}]$, 
  which in our notations is  $Y_{\wK''}$. 
  By  \cite[\S 7]{kwlan-book}, $[S_{\wK''}]$    is an algebraic space and even a quasi-projective scheme. Moreover, by a general property of moduli stacks (see \cite[Chap.~7]{olson-book}), $\phi$ is an isomorphism over the locus $U$ where there is {\it no non-trivial automorphism}, by which we mean it has no infinitesimal automorphism; analytically, this corresponds to  points of $Y_{\wK}(\C)$ having no non-trivial stabilizers. 
Now $U$ is a priori an open subscheme of $[S_{\wK''}]$, but since it is where $\phi$ is an isomorphism, we get a canonical open  $U \hookrightarrow S_{\wK''}$. This  map  tautologically yields the desired family $f: A\to U$ of abelian varieties of Picard type. 
In view of Lemma~\ref{lem:reflexions} and the discussion preceding it, after possibly enlarging $k$, we may consider $U$ be the complementary open in $Y''$ of its finitely many elliptic points. 
\end{proof} 

\vspace{-8mm}

\subsection{\'Etale fundamental groups and Mumford--Tate groups}\label{MT-groups}
Let $k$ be a number field  containing $M$ over which the connected component of $Y_{\wK}$ are defined.
Fix a connected component $Y$ of $Y_{\wK}\times_M k$ and a smooth open  $U$ of $Y$ endowed with an abelian scheme $f:A\to U$ of Picard type. Denote by  $\eta$ the generic point of the smooth surface $U$.
Fixing a closed geometric point $\bar x$ of $U$ the \'etale fundamental group  sits in the middle of a short exact sequence
\begin{align}\label{fund-group-es}
1\to \Pi_1(U_{\bar k}, \bar x)\to \Pi_1(U, \bar x) \to \Gal_k \to 1.
\end{align}
The morphism  $f:A\to U$  being proper and smooth, one can consider the  \'etale sheaf $\mathrm{R}^1f_*\Z_\ell$ on $U$.
As $U$ is geometrically connected, we have $\Pi_1(U, \bar x)\simeq  \Pi_1(U, \bar \eta)$ and the latter acts on
\[(\mathrm{R}^1f_*\Z_\ell)_{\bar \eta}= \rH^1(A_{\bar \eta}, \Z_\ell)=(T_\ell A_{\eta})^\vee,\]
yielding a continuous representation
\[\Gal(\bar \eta/\eta)\twoheadrightarrow \Pi_1(U, \bar x)\xrightarrow{\rho_{U,\ell}} \Aut_{\Z_\ell}(T_\ell A_{\eta}).\]
Any closed point $x\in U(k)$ yields a section $s_x:\Gal_k\to  \Pi_1(U, \bar x)$ of \eqref{fund-group-es} allowing one to consider
\[\rho_{x,\ell}=\rho_{U,\ell}\circ s_x: \Gal_k \to \Aut_{\Z_\ell}(T_\ell A_{\eta}).\]
Finally for any closed curve $C\subset U$ defined over $k$, there is a natural map   $\Pi_1(C, \bar x)\to  \Pi_1(U, \bar x)$
whose composition with  $\rho_{U,\ell}$ is denoted $\rho_{C,\ell}$. 
As  $f:A\to U$ is of Picard type, for any $x\in C(k)$,
\[\Gamma_x=\mathrm{im}(\rho_{x,\ell})\subset\Gamma_C=\mathrm{im}(\rho_{C,\ell})\subset \Gamma_U=\mathrm{im}(\rho_{U,\ell})\subset K^\circ_\ell.\]

By a series of results of Cadoret--Tamagawa
(see \cite{cadoret-tamagawa-inv, cadoret-tamagawa-duke1}), the set $C_\rho$ of all   $x\in C(k)$ for which $\Gamma_x$ is not open in $ \Gamma_C$ is finite and the index $[\Gamma_C: \Gamma_x]$ is uniformly bounded for  $x\in C(k){\smallsetminus} C_\rho$.

Recall that the Mumford--Tate group $\MT(A)$ of a polarized abelian variety $A$ over $\C$  is 
the smallest connected reductive subgroup of $\GL(\rH_1(A, \Q))$ over $\Q$, 
whose $\R$-points contain the image of the associated $\R$-morphism $h: \C^\times \to \GL(\rH_1(A(\C), \R))$ coming from the Hodge decomposition. 
If  $A$ is defined over a  field  $k$  finitely generated  over $\Q$, then the image $\Gamma_\ell$ of $\Gal_k$ acting on $T_\ell A$ is an $\ell$-adic Lie group. 
Denoting  $\mathfrak{g}^{}_{\Z_\ell}$ its Lie algebra, 
it is a theorem of Deligne  \cite[Chap. I.2]{deligne-milne-ogus-shih} that 
$\mathfrak{g}^{}_{\Q_\ell}=\mathfrak{g}^{}_{\Z_\ell}\otimes_{\Z_\ell}\Q_\ell \subset \Lie(\MT(A)\otimes_{\Q}\Q_\ell)$ 
and the Mumford--Tate conjecture, known for
abelian varieties of dimension at most $3$,  asserts that they are equal (see {\it e.g.} \cite{chi}).

As the  Mumford--Tate group of the (generic point of the) universal family $f:A\to U$  is given by 
$\wG=\mathrm{GU}(2,1)$, it follows from Deligne \cite[Cor.~6.2]{deligne-milne-ogus-shih} that 
  the Mumford--Tate group of any abelian $3$-fold  of Picard type is a reductive subgroup of   $\wG$.
We have the following trichotomy. 

\break 
\begin{lemma}\label{prop:tricotomy}
Let $\fg$ be a reductive Lie subalgebra of $\mathfrak{gu}(3,\Q_\ell)$. 
If $\fg'\subset \mathfrak{su}(3,\Q_\ell)$ denotes the semi-simple   part   of $\fg$,   exactly one of the following holds:
\begin{enumerate}[wide]
\item $\fg'=\{0\}$, {\it i.e.}, $\fg$ is abelian, 
\item  $\fg'$ is a form of $\mathfrak{sl}(2,\Q_\ell)$, 
\item  $\fg'=\mathfrak{su}(3,\Q_\ell)$. 
\end{enumerate}
\end{lemma}

\begin{proof}
If $\fg'_{\overline\Q_\ell}=\{0\}$, then $\fg'=\{0\}$, whereas if $\fg'_{\overline\Q_\ell}=\mathfrak{sl}(3,\overline\Q_\ell)$,  then $\fg'=\mathfrak{su}(3,\Q_\ell)$ for dimension reasons.  In the remaining cases, using the well known fact that any proper non-zero semi-simple Lie subalgebra of $\mathfrak{sl}(3,\overline\Q_\ell)$ is isomorphic to $\mathfrak{sl}(2,\overline\Q_\ell)$, we deduce that $\fg'$ is a form of $\mathfrak{sl}(2,\Q_\ell)$.
\end{proof}

If $A$ is of CM type (resp. admits  a   CM factor), then $\Lie(\MT(A)\otimes_{\Q}\Q_\ell)$ is of the first (resp. second) type in Lemma~\ref{prop:tricotomy}. 
Also clearly  $\rho_{A,\ell}$ is potentially abelian (resp. is potentially  reducible)
if, and only if,  $\mathfrak{g}^{}_{\Q_\ell}$ is of the first (resp. second) type in Lemma~\ref{prop:tricotomy}. 
The following proposition bridges the two sides.

\begin{prop}\label{prop:MT}
Let $A$ be an abelian $3$-fold  of Picard type over a number field $k$. 
  Then $\rho_{A,\ell}$  is potentially abelian (resp. is potentially  reducible)  
  if, and only if, $A$ is  of  CM type (resp. admits  a  CM factor).
\end{prop}

\begin{proof} One implication is clear. For the other, after extending scalars to $\overline\Q_\ell$ 
and after  possibly replacing $k$ by a finite extension, we may assume that 
 $\rho_{A,\ell}$ is potentially reducible, {\it i.e.} it contains (as a direct factor by Faltings)
a character  $\chi^{}_\ell: k^\times\backslash \A_k^\times \to \overline\Q_\ell^\times$. Being a sub-representation of  $\rho_{A,\ell}$, $\chi^{}_\ell$ is unramified outside a finite set of places, 
and its restriction to decomposition groups  at places above $\ell$ is Hodge-Tate with weights $0$ and $-1$. 
In addition it is pure of weight $-1$. By Weil,  $\chi^{}_\ell$ corresponds to an algebraic Hecke character  $\chi: k^\times\backslash \A_k^\times \to \C^\times$ whose infinity component is necessarily of the form
$\rN_{\Phi'}\circ \rN_{k/L'}$ where $\rN_{\Phi'}$ is the partial norm given by a CM type $\Phi'$ for a CM field
$L'\subset k$. By \cite[Lem.~2]{shimura-CM} replacing $(L',\Phi')$ by its double reflex yields the same infinite component, hence we may and do assume that $(L',\Phi')$ is  primitive, {\it i.e.}, coincides with  the reflex of a CM field $L$
endowed with a CM type $\Phi$. Further replacing $k$ by a finite (abelian) extension one can assume that
$\chi^{}_f$ takes values in  $L^\times$.
 By  Casselman (see \cite[Thm.~6]{shimura-CM}), there exists an  abelian variety $A'$
  defined over $k$  which is CM of type $(L,\Phi)$ and such that  $\rho_{A',\ell}=\chi^{}_\ell$, hence
\[\Hom_{\Gal_k}(\rho_{A,\ell}, \rho_{A',\ell})\neq \{0\}.\]
By Faltings, one deduces that $\Hom_{k}(A,A')\ne  \{0\}$, hence $A$ contains a non-trivial  CM quotient. 
If $A$ is of CM type  we're done. If not, then  there exists  an abelian variety $A''$ which is not of CM type and such that  $A$ is isogenous to $A'\times A''$, {\it i.e.}, $V_{\ell}A=V_{\ell}A'\oplus V_{\ell}A''$. Furthermore since $\Hom_{k}(A',A'')=\{0\}$, one can show that the  isogeny is $\cO_M$-linear. Hence $A''$
is also of Picard type,   from which one  deduces that we're in the second case in Lemma~\ref{prop:tricotomy}, in particular $\rho_{A,\ell}$  is not potentially abelian. 
\end{proof}

\section{Uniform irreducibility of Galois images}\label{s:uniform}

\subsection{A  lemma on surfaces with isolated singularities}
We denote by $q(X)$ the irregularity of a projective  algebraic surface
$X$,  given by the dimension of $\rH^1(X, \cO_{X})$.
If $X$ is smooth and projective over $\C$, then $q(X)=\dim \rH^0(X, \Omega^1_{X})$.

Let $X$ be a projective irreducible algebraic surface over $\C$ with isolated singularities, {\it i.e.}, such that there exists a
smooth  open $j:U\hookrightarrow X$ whose complement   $Z=X{\smallsetminus} U$ consists of finitely many closed points.
There exists a smooth resolution
\[\phi: \widetilde{X}  \to X\]
such that $\phi^{-1}(Z)$ is a divisor with normal crossings with $\phi$ restricting to   an isomorphism from
 $\phi^{-1}(U)$ onto $U$.   Thus we get an injection
$\widetilde j: U  \hookrightarrow  \widetilde{X} $ such that $j= \phi \circ\widetilde j $, and we denote by
\[\widetilde j^\ast: \rH^1(\widetilde{X}, \Q)  \to  \rH^1(U, \Q)\]
the pullback homomorphism on Betti cohomology. By \cite[Thm.~3.2.5(iii)]{deligne-hodgeII} we know that  $\widetilde j^\ast$ is a homomorphism of
mixed Hodge structure, with $ \rH^1(\widetilde{X}, \Q)$ being   pure of weight $1$.

\begin{lemma} \label{top} The map $\widetilde j^\ast$ is an isomorphism, in particular  $\rH^1(U, \Q)$ is a pure weight $1$ Hodge structure
and $q(X)=\dim \rH^0(U, \Omega^1_U)$.
\end{lemma}

\begin{proof}
Let $\mathrm{IH}^\bullet(X,\Q)$ denote the middle intersection cohomology of $X$. Since
\[\dim(X)-1>0=\dim(Z)\]
by  \cite[Thm.~5.4.12]{dimca}
$j^\ast: \mathrm{IH}^1( X, \Q)  \to  \mathrm{IH}^1(U, \Q)$ is an isomorphism, while $\widetilde j^\ast$
is injective. Moreover by Corollary~5.4.11 and Proposition~5.4.4 in {\it loc.cit.}
$\phi^\ast: \mathrm{IH}^1(X, \Q)  \to  \mathrm{IH}^1(\widetilde{X}, \Q)=\rH^1(\widetilde{X}, \Q)$  is an embedding, while $\mathrm{IH}^1( U, \Q)=\rH^1( U, \Q)$.
This is summarized in the following commutative diagram:
\[\xymatrix@C=40pt{
\mathrm{IH}^1( X, \Q) \ar@{^{(}->}_{\phi^\ast}[d]   \ar_{\sim}^{j^\ast}[r] &\mathrm{IH}^1( U, \Q) \ar@{=}[d] \\
 \rH^1( \widetilde{X}, \Q) \ar@{^{(}->}^{\widetilde j^\ast}[r]   & \rH^1( U, \Q) }\]
It immediately follows that $\widetilde j^\ast$ is an isomorphism and $q(\widetilde{X})=
\dim \rH^0(\widetilde{X}, \Omega^1_{\widetilde{X}})= \dim \rH^0(U, \Omega^1_U)$. Finally $q(X)=q(\widetilde{X})$ as the
irregularity is a birational invariant. \end{proof}

\subsection{Irregularity of Picard modular surfaces}\label{irregularity}
Let  $z\mapsto \bar z$ be the non-trivial automorphism of $M/\Q$.
Put  $M^1=\{z\in M^\times\mid z\bar z=1\}$ viewed as an algebraic torus over $\Q$ and denote by  $\A_M^1$ its adelic points.

For $\wK\subset \wG(\A_f)$ an open compact subgroup, we recall  the Shimura variety of Picard type $Y_{\wK}$ from 
\S\ref{sec:moduli}.  Let $G^1=\mathrm{SU}(V)$ be the derived group of $\wG$.  As $G^1$ is simply connected and $G^1_\infty$ is not compact, the  Strong Approximation Theorem (see~\cite[Thm.~7.12]{platonov-rapinchuk}) implies that
$G^1(\Q)$ is dense in $G^1(\A_f)$. It follows that the
determinant map defines an isomorphism between the group of connected components  $\pi_0(Y_{\wK})$  and the  idele class group
$\A_M^\times/M^\times \det(\wK)M_\infty^\times$. Shimura's theory of canonical models implies that the connected components of
$Y_{\wK}$ are all Galois conjugates, hence share the same irregularity, and the same is true for the Shimura variety
for $G$
\[Y_K(\C)= G(\Q)\backslash \left(\mathcal{H} \times G(\A_f)/K\right), \]
where $K=\wK\cap G(\A_f)$. Letting $\Gamma= \wG(\Q)\cap  \wK \wG(\R)$
 it follows from $\nu(\Gamma)\subset \Q^\times \cap  \widehat\Z^\times \widehat\R_+^\times=\{ 1\}$  that both
$Y_{\wK}$ and $Y_K$ share the same connected component of identity given by $Y_\Gamma=\Gamma\backslash \mathcal{H}$
(see \cite[(8)]{dimitrov-ramakrishnan}). One should be careful to observe that the natural dominant map $Y_{K^1}\to Y_\Gamma$, where $Y_{K^1}$  is
the  Shimura variety of level $K^1=K\cap G^1(\A_f)$ for  $G^1$, is an isomorphism precisely when either $\det(\Gamma)=\{1\}$ or $-1\in \Gamma$.

\begin{prop} The irregularity of any connected component of the minimal compactification $Y_{\wK}^\ast$ of $Y_{\wK}$ is given by the formula
\begin{equation}\label{q-formula}
q(Y_\Gamma^\ast)= \sum_{(\lambda,\nu)\in \Xi/\widehat{\pi_0(Y_K)}}
\sum_{ \pi_f \in \Pi_f(\lambda,\nu)} \dim(\pi_f^K) \frac{1+W(\lambda\nu^{-1}_M)(-1)^{s(\pi_f)}}{2}, \text{ where }
 \end{equation}
\begin{itemize}
\item  $\Xi$  is the set of  pairs $(\lambda,\nu)$ of a   Hecke character $\lambda$ of $M$ whose restriction to $\Q$ is $\left(\tfrac{\cdot}{D}\right)$ and of a   character  $\nu$   of  $\A_M^1/M^1$, such that
\begin{equation*}
\lambda_\infty(z)= \frac{\bar{z}}{|z|}\text { , for all } z\in M_\infty^\times\simeq\C^\times,
\text { and }  \nu_\infty(z)= z,\text{ for all } z\in M_\infty^1,
\end{equation*}
\item $\Pi_f(\lambda,\nu)$ is the finite part of the global Arthur packet from  \S\ref{packets},
\item  $W(\lambda\nu_M^{-1})\in \{\pm 1\}$  is the global root number, where  $\nu^{}_M(z)= \nu(z/ \bar z)$ for $z\in \A_M^\times$,
\item $s(\pi_f)$  the number of finite  places $v$  at which  $\pi_v\simeq \pi_{c}(\lambda_v,\nu_v)$,  and
\item $\mu \in \widehat{\pi_0(Y_K)}$ acts freely on $\Xi$ by sending $(\lambda,\nu)$ to $(\lambda\mu^{}_M,\nu\mu)$.
\end{itemize}
\end{prop}
\begin{proof} We first observe that the complement of $Y_\Gamma$ in $Y_\Gamma^\ast$ consists of finitely many points,  the cusps. 
Furthermore, as explained in \S\ref{sec:abelian-family},  $Y_\Gamma$ admits  finitely many isolated singularities all of which are  elliptic points. 
Thus  there exists a smooth open $U_K$ of the normal projective surface $Y_K^\ast$  whose complement  consists of finitely many closed points.
Lemma~\ref{top} applied  component-wise  to $U_K$ yields $q(Y_K^\ast)= \dim\rH^0(U_K, \Omega^1_{U_K})$.
Let   $K'$ be any normal finite index torsion free subgroup of $K$, {\it e.g.}
the intersection with the principal congruence subgroup of level $3$ (see  \cite[Lem.~1.4]{dimitrov-ramakrishnan}). By Koecher's Principle, as
$Y_{K'}{\smallsetminus} U_{K'}$ has codimension at least  $2$  in $Y_{K'}$, we have
\[\dim\rH^0(U_K, \Omega^1_{U_K})=
\dim\rH^0(U_{K'}, \Omega^1_{U_{K'}})^{K/K'}=\dim\rH^0(Y_{K'}, \Omega^1_{Y_{K'}})^{K/K'},\]
where $U_{K'}$ is the inverse image of $U_{K}$ under the natural projection $Y_{K'}\to Y_{K}$.
Taking invariants by the finite group $K/K'$ in Rogawski's formula \cite[(14)]{dimitrov-ramakrishnan} for $q(Y_{K'}^\ast)=\dim\rH^0(Y_{K'}, \Omega^1_{Y_{K'}})$ yields
\[\dim\rH^0(U_K, \Omega^1_{U_K})= \sum_{(\lambda,\nu)\in \Xi}
\sum_{ \pi_f \in \Pi(\lambda_f,\nu_f)} \dim(\pi_f^K) \frac{1+W(\lambda\nu_M^{-1})(-1)^{s(\pi_f)}}{2}\]
(there is a misprint in  {\it loc. cit.} where one should read $(1+W(\lambda\nu_M^{-1})(-1)^{s(\pi_f)})$ instead of
$(W(\lambda\nu_M^{-1})+(-1)^{s(\pi_f)})$; it also uses the inverse notation for the base change $\nu^{}_M$). The presence of this root number translates the fact that for $\pi_f \in \Pi(\lambda_f,\nu_f)$ and
$\pi_\infty$ the unique non-tempered holomorphic  representation in the local Arthur packet $\Pi(\lambda_\infty,\nu_\infty)$,
$\pi=\pi_f \otimes \pi_\infty$ is automorphic if, and only if,  $W(\lambda\nu_M^{-1})=(-1)^{s(\pi_f)}$.
Both  $\dim(\pi_f^K)$ and $1+W(\lambda\nu_M^{-1})(-1)^{s(\pi_f)}$ being preserved by the action of  $\widehat{\pi_0(Y_K)}$,
one deduces the desired formula for $q(Y_\Gamma^\ast)$ as in \cite[(15)]{dimitrov-ramakrishnan}.
\end{proof}

\subsection{Twists of canonical characters and root numbers}\label{canonical}
Hecke characters  $(\lambda,\nu)\in \Xi$ whose local components at each finite place have `minimal' ramification are intimately related to the canonical characters studied by Gross and Rohrlich. They  play a pivotal role in our production of automorphic forms contributing to the irregularity of the Picard modular surfaces of low level. We will now briefly recall some of their properties under the running assumption that  $D> 3$ is odd. Consider the character $\lambda_\infty(z)=\bar{z}\cdot  |z|^{-1}$ of $M_\infty^\times\simeq \C^\times$
and let $\lambda_f:\widehat{\cO}_{M}^\times\to \C^\times$ be  a continuous  character  whose restriction  to $\cO_{M,p}^\times$
is given by the unique quadratic character
\[\cO_{M,p}^\times\to \F_p^\times\xrightarrow{\left(\tfrac{\cdot}{p}\right)}\{\pm 1\},\]
 for all $p$ dividing  $D$, and is trivial otherwise. As $\left(\tfrac{-1}{D}\right)=-1$, it follows that $\lambda_\infty$ and $\lambda_f$ agree on  $\cO_M^\times=\{\pm 1\}$.
The finiteness of the ideal class group $\Cl_M=\A_M^\times/M^\times \widehat{\cO}_{M}^\times M_\infty^\times$ guarantees that the resulting character of  $M^\times \widehat{\cO}_{M}^\times M_\infty^\times$ extends to a character $\lambda$ of $\A_M^\times$ and clearly two such extensions must differ by an ideal  class character. As by construction the restriction of $\lambda_f$ to $\widehat{\Z}^\times$ agrees with
the  quadratic Dirichlet character  $\left(\tfrac{\cdot}{D}\right)$ viewed as a character of $\A^\times/\Q^\times \mathrm{Nm}(\A_M^\times)=\Gal(M/\Q)$ and  $\A^\times = \Q^\times  \widehat{\Z}^\times \R_+^\times$,
 it follows that the restriction of $\lambda$ to  $\A^\times$ equals  $\left(\tfrac{\cdot}{D}\right)$, 
 {\it i.e.}, $\lambda$ is conjugate-symplectic.
 Such a character is called  canonical  and we will denote it by $\lambda_c$, remembering that it is only unique up to a multiplication by a  character of  $\Cl_M$.   The root number $W(\lambda_c^3)=-W(\lambda_c)=\left(\tfrac{-2}{D}\right)$ is $1$  if, and only if, $D\equiv 3\pmod{8}$ (see  \cite{rohrlich-canonical}).

Assume henceforth that $\det(K)=\widehat{\cO}_{M}^1$, so that $\pi_0(Y_K)=\Cl_M^1:=\A_M^1/M^1 \widehat{\cO}_{M}^1 M_\infty^1$, and  that $\widehat{\cO}_{M}^1$ embeds centrally into $K$, the central character $\omega=\nu\cdot \lambda_{|M^1}$ of any
$\pi$ contributing to \eqref{q-formula} has to be everywhere unramified, {\it i.e.},
\begin{equation}\label{q-formula-can}
q(Y_\Gamma^\ast)= \frac{1}{|\Cl_M^1|}
\sum_{\chi\in \Xi^1, \omega\in \widehat{\Cl_M^1}}
\sum_{ \pi_f \in \Pi_f\left(\lambda_{c}\chi^{}_{M},{\lambda^{-1}_{c}}_{\!\!\!\!\!|M^1} \chi^{-2}\omega\right)} \dim(\pi_f^K) \frac{1+W(\lambda_c^3\chi^{3}_{M})(-1)^{s(\pi_f)}}{2},
\end{equation}
where $\Xi^1$ denotes the set of finite order characters of $\A_M^1/M^1$ (see \cite[(3)]{rohrlich-canonical} for the fact that multiplication by an ideal class character does not change the root number).

If $3$ does not divide the class number  $h:=|\Cl_M|$, then the action of  $\mu \in \widehat{\Cl_M^1}$ sending $(\chi,\omega)$ to
$(\chi\mu,\omega\mu^3)$ allows one to twist out the central character $\omega$ and obtain the simpler formula
\begin{equation}
q(Y_\Gamma^\ast)= \sum_{\chi\in  \Xi^1}
\sum_{ \pi_f \in \Pi_f(\lambda_{c}\chi^{}_{M})} \dim(\pi_f^K) \frac{1+W(\lambda_c^3\chi^{3}_{M})(-1)^{s(\pi_f)}}{2}.
\end{equation}

Successfully applying  \eqref{q-formula-can} requires one to understand how root numbers behave under twisting.
As we are interested in creating irregularity at  level $\Gamma''_0(\ell^r)$, we focus on characters $\chi$ which are 
only ramified at the fixed  inert  prime $\ell$.

\begin{lemma}\label{lem:twisting} 
For  $\chi\in \Xi^1$ of Artin conductor $\ell^{a(\chi)}$ and $\lambda$ an odd power of a canonical character,
\[ W(\lambda\chi^{}_{M}) =(-1)^{a(\chi)}  \chi^{}_\ell(-1) W(\lambda).\]
\end{lemma}

\begin{proof} Using the factorization of root numbers $W = \prod\limits_v W_v$, it suffices to prove that
\begin{align} 
\label{root-at-ell}  W_\ell(\lambda_\ell\chi^{}_{M,\ell}) &  = (-1)^{a(\chi)}  \chi^{}_\ell(-1) W_\ell(\lambda_\ell) \text{ and } \\
\label{root-not-ell} W_v(\lambda_v\chi^{}_{M,v}) & = W_v(\lambda_v), \, \text{ for all  }  v \ne \ell, 
\end{align}
where the local  factors are defined  using  the standard additive character $\psi^{}_M=\psi^{}_\Q\circ \Tr_{M/\Q}$. Applying  \cite[Prop.~3]{rohrlich-galois}  to both $\lambda_\ell$ and $\lambda_\ell\chi^{}_{M,\ell}$ yields  \eqref{root-at-ell}.
As $\chi^{}_{M, \infty} =1$, it suffices to check \eqref{root-not-ell} for $v$   finite. Moreover, the characters  $\lambda_v$ and $\chi^{}_{M,v}$ are unramified for  $v\nmid\ell D$, hence  both sides of \eqref{root-not-ell}  are $1$. 
Finally, for  $v$ dividing $D$,  $\chi^{}_{M,v}$ is unramified while   $\lambda_v$ is tamely ramified and $\psi^{}_{M_v}$ has conductor $1$, implying by \cite[(5.5.1)]{deligne-epsilon} that
$W_v(\lambda_v\chi^{}_{M,v})/ W_v(\lambda_v)=\chi^{}_{M,v}(\xi^{1+1\cdot 1})=1$.
\end{proof}

\subsection{The Bombieri--Lang Conjecture for Picard modular surfaces}
In this part we follow the general strategy of \cite{dimitrov-ramakrishnan} by adapting it to the case where the level is not neat, the main point being to show that the irregularity of the Picard modular surfaces under consideration is at least $3$.
This requires some new techniques,  also providing new proofs to the cases treated in {\it loc.~cit.}

We recall the  index $2$ subgroup $\wK''$ of  the maximal open compact subgroup   $\wK^\circ $ of $\wG(\A_f)$ introduced in \eqref{K-double-prime}. Given  $r\in\Z_{\geqslant 1}$ and a prime $\ell$ inert in $M$,  we let $\wK''_0(\ell^r)$  be the subgroup of $\wK''$  whose component at $\ell$ is the depth $r$ Iwahori subgroup $I_r$  and we let   
  \[K''_0(\ell^r)=\wK''_0(\ell^r)\cap G(\A_f) \text{  and }   \Gamma''_0(\ell^r)= G(\Q) \cap K''_0(\ell^r)\cdot G(\R).\]

Consider an  automorphic representation $\pi\in\Pi\left(\lambda_{c}\chi^{}_{M},{\lambda^{-1}_{c}}_{\!\!\!\!\!|M^1} \chi^{-2}\omega\right)$ having non-zero  $K''_0(\ell^r)$-invariants. By Propositions~\ref{ram-inv} and \ref{prop:inert-inv}, for all finites place $v\ne \ell$ we must have  $\pi_v=\pi_{n,v}$ with  $\chi^{}_v$ unramified, in which case  it follows from   \eqref{KD-double-prime} that  the line $\bigotimes\limits_{p\mid D} \pi_{n,p}^{K''_p}$ is  fixed by $K''_D$.  By  \eqref{q-formula-can},   $\pi$ contributes to $q\left(Y_{\Gamma''_0(\ell^r)}^\ast\right)$ 
if either   $\pi_\ell=\pi_{n,\ell}$ and $W(\lambda_c^3\chi^{3}_{M})=1$, or  $\pi_\ell=\pi_{c,\ell}$ and $W(\lambda_c^3\chi^{3}_{M})=-1$. 
As  the choice of a  character  $\chi^{}_\ell:\Q_{\ell^2}^1 \to \C^\times$  determines uniquely, up to an ideal class character, a global character  $\chi^{}_{M}$ unramified outside $\ell$, the irregularity formula \eqref{q-formula-can} becomes (we have  omitted $\omega\in \widehat{\Cl_M^1}$ as  it is trivial on $\Q_{\ell^2}^1$)
 \begin{equation}\label{q-formula-at-ell}
\tfrac{1}{h}\cdot q\left(Y_{\Gamma''_0(\ell^r)}^\ast\right)=  
\sum_{\chi^{}_\ell, \text{ s.t.} W(\lambda_c^3\chi_M^3)=1} \dim(\pi_{n,\chi^{}_\ell}^{I_r}) + 
\sum_{\chi^{}_\ell, \text{ s.t.} W(\lambda_c^3\chi_M^3)=-1} \dim(\pi_{c,\chi^{}_\ell}^{I_r}),  
\end{equation}
where $\Pi(\lambda_0 \chi^{}_{M,\ell}, \chi^{-2}_\ell)=\{\pi_{n,\chi^{}_\ell},\pi_{c,\chi^{}_\ell}\}$ denotes the local Arthur packet attached to the  quadratic unramified character $\lambda_0:\Q_{\ell^2}^\times \to \C^\times$ and to $\chi^{}_\ell:\Q_{\ell^2}^1 \to \C^\times$.

\begin{rem} \label{odds-ends} 
The study of the conductor of the supercuspidal non-generic representation $\pi_{c,\chi^{}_\ell}$ appears to be very delicate.  Indeed, even in the depth $0$  case ({\it i.e.}, trivial $\chi$), some preliminary computations suggest that it does not contain non-zero $K_T$-invariant vectors. 
In a previous paper (\cite{dimitrov-ramakrishnan}),  we mistakenly switched $\pi_{c}$  and $\pi_{2}$ (which are in the same $L$-packet)  in the  proof of Propositions~3.6 and 3.8 in {\it loc. cit.}. It did not affect Theorem~0.2 there at all, but in Theorem~0.1 when $W(\lambda^3)=(-1)^d$, 
we need to replace $\Gamma_1(\mathfrak{C})$ by $\Gamma_1(\mathfrak{Cq})$ for any   $\mathfrak{q}$ prime of $F$ relatively prime to 
 $\mathfrak{C}$ and non-split in $M$. We overcome the difficulty caused by $\pi_{\mathfrak{q},c}$ not having $\Gamma_1(\mathfrak{q})$-invariants by using an appropriate  twist $\lambda_c\chi^{}_{M}$ of $\lambda_c$ for which  $W(\lambda_c^3\chi^{3}_{M})=(-1)^{d-1}$ and the corresponding  $\pi_{\mathfrak{q},n}$ has the requisite invariants. Adapting this twisting method at ramified places, allows us  to establish a variant of  Theorem~0.3 from  {\it loc. cit.}. 
This will be taken up elsewhere. 
\end{rem}

Henceforth, we will only use  a lower bound for the irregularity corresponding to the contribution of  everywhere non-tempered automorphic representations. Namely, 
combining \eqref{q-formula-at-ell} with Lemma~\ref{lem:twisting}, we have 
 \begin{equation}\label{q-inequality-at-ell}
 q\left(Y_{\Gamma''_0(\ell^r)}^\ast\right)\geqslant h \cdot  \sum_{\chi^{}_\ell} \dim(\pi_{n,\chi^{}_\ell}^{I_r}), 
\end{equation}
where the sums runs over characters  $\chi^{}_\ell:\Q_{\ell^2}^1 \to \C^\times$ such that 
$\chi^{}_\ell(-1)=(-1)^{a(\chi^{}_\ell)+(D-3)/4}$.

\begin{prop} \label{prop:q3} We recall that  $D>3$ is odd and   $h$ denotes the class number of $M=\Q(\sqrt{-D})$.
\begin{enumerate}[wide]
\item If $D\equiv 3\pmod{8}$ then $q\left(Y^\ast_{\Gamma_0''(\ell)}\right)=q(Y^\ast_{\Gamma''})=h$ and
$q\left(Y^\ast_{\Gamma_0''(\ell^{2r})}\right)\geqslant (r+1)h$, for $r\in \Z_{\geqslant 1}$. 

\item If $D\equiv 7\pmod{8}$ then $q\left(Y^\ast_{\Gamma_0''(27)}\right)\geqslant h$, 
 $q\left(Y^\ast_{\Gamma_0''(5^6)}\right)\geqslant  3h $ and 
$q\left(Y^\ast_{\Gamma_0''(\ell^3)}\right)\geqslant  3h $,   for $\ell\geqslant 7$.   
\end{enumerate}
\end{prop}

\begin{proof}  
(i) We can take $\chi=\mathbf{1}$ in  \eqref{q-inequality-at-ell}. The claim follows from 
Proposition~\ref{prop:inert-inv} and Corollary~\ref{cor:deep-iwahori-inv}. 

(ii) In this case, we use tamely ramified $\chi_\ell$  in \eqref{q-inequality-at-ell}, {\it i.e.} $a(\chi^{}_\ell)=1$,  such that $\chi^{}_\ell(-1)=1$. 
For $\ell\geqslant 3$, there are precisely $\frac{\ell-1}{2}$ such characters, if $3 \nmid (\ell+1)$,  and $\frac{\ell-5}{2}$ choices,  if $3 \mid (\ell+1)$.  In particular, there are at least $3$ choices for all  $\ell\geqslant 7$.  
As $K_T$ contains a conjugate $I_3$, the claim then follows from Theorem~\ref{inert-inv-ram} 
supplemented, when  $\ell=5$, by a  numerical computation showing that $\pi_{n,\chi^{}_5}$ has  non-zero $I_6$-invariants for $\chi^{}_5$ such that $\chi^{}_5(-1)=-1$ and 
 $a(\chi^{}_5)=2$.  
 \end{proof}

The following results  from Faltings \cite{faltings-lang}  (see \cite[\S 3.2]{dimitrov-ramakrishnan}  for details).

\begin{thm}  Let $\wK$ an open compact subgroup   of $\wG(\A_f)$ and let $\Gamma= \wG(\Q)\cap  \wK \wG(\R)$.
If $q\left(Y^\ast_{\Gamma}\right)\geqslant 3$, then  $Y^\ast_{\wK}$ satisfies the   Bombieri--Lang Conjecture.
\end{thm}

As $2$ splits in $M$ for $D\equiv 7\pmod{8}$, Proposition~\ref{prop:q3} has the following consequence.

\begin{cor} \label{bomb-lang}
If $D\equiv 3\pmod{8}$,   then the Bombieri--Lang Conjecture holds for $Y^\ast_{\wK''(\ell^4)}$, and even for
$Y^\ast_{\wK''}$, when $h\geqslant 3$.
If $D\equiv 7\pmod{8}$,  then the  Bombieri--Lang Conjecture holds for 
$Y^\ast_{\wK''(\ell^3)}$ for $\ell\geqslant 7$, and also for $Y^\ast_{\wK''(3^7)}$ and  $Y^\ast_{\wK''(5^{6})}$. 
\end{cor}

\subsection{Proof of Main Theorem}
At  different stages of the proof we will remove finite sets and deal with them in the last step.
By  \S\ref{sec:moduli}  an abelian variety $A$ as in the Theorem~\ref{theoremB} defines a $k$-rational point on $Y_{\wK''}$, where $\wK''$ is
defined in \eqref{K-double-prime}.
If $A[\ell^r]$ admits a full $k$-rational flag (or equivalently, a $k$-rational isotropic line),  we claim that $A$ defines a $k$-rational point on $Y_{\wK''_0(\ell^r)}$,
where $\wK''_0(\ell^r)\subset\wK''$ is the  subgroup whose component at $\ell$ consists of elements whose reduction  modulo $\ell^r$ belong to  the standard Borel subgroup  $\widetilde{B}(\Z/\ell^r\Z)$ of $\wG(\Z/\ell^r\Z)$. Indeed,  by assumption one  knows that there is
 some Borel subgroup containing  the image of $\Gal_k$ acting on $A[\ell^r]$. However $\widetilde{G}(\Q_\ell)$ acts transitively on isotropic lines
(because isometries between hermitian subspaces always extend), hence all Borel subgroups are conjugated by $\widetilde{G}(\Q_\ell)$, and in fact by  $\wK_\ell^\circ=\widetilde{G}(\Z_\ell)$ (using Iwasawa  decomposition).  As  $\wK''$ is a normal subgroup of $\wK^\circ$ we deduce that  the Galois image is contained in $\wK''\cap \wK_0(\ell^r)=\wK''_0(\ell^r)$, proving the claim.

 By   Corollary~\ref{bomb-lang}, $Y^{\ast}_{\wK''_0(\ell^{7})}$  satisfies the Strong Bombieri--Lang Conjecture.
 In particular all  its $k$-rational points lie in a subvariety $Z$ defined over $k$ which is a finite union of points and curves.
 
Let us now take one of the (finitely many) geometrically connected curve  $C$ in $Z$, and after removing finitely many of its points
(which would not affect the wanted result), we may assume that  $C$ is contained in the smooth open $U$ from \S\ref{sec:abelian-family}. In particular, there exists a family $f:A\to C$ of abelian $3$-folds of Picard type.

As in \S\ref{MT-groups} let   $\Gamma_C\subset \wK_\ell^\circ=\widetilde{G}(\Z_\ell)$ be the image of the \'etale fundamental group acting on the $\ell$-adic Tate module of the generic fiber of the family. By Cartan's theorem (see \cite[LG5.42]{serre-LALG}),  $\Gamma_C$ is an $\ell$-adic Lie group hence admits a  Lie algebra $\mathfrak{g}^{}_{\Z_\ell}$. By Bogomolov \cite{bogomolov-SSSR} the Lie algebra $\mathfrak{g}^{}_{\Q_\ell} =\mathfrak{g}^{}_{\Z_\ell}\otimes_{\Z_\ell}\Q_\ell$ is algebraic, namely it is the Lie algebra of the Zariski closure of $\Gamma_C$ in $\widetilde{G}(\Q_\ell)$, the latter being   furthermore  reductive over $\Q_\ell$ by   Faltings \cite[Thm.~3]{faltings-AG}.
By the Mumford--Tate Conjecture, which is known for is known for abelian $3$-folds (see {\it e.g.} \cite{chi}), we know that
$\mathfrak{g}^{}_{\Q_\ell}$ is the Lie algebra of the Mumford--Tate group $\MT(A)\otimes_{\Q} \Q_\ell$.
 As $C$  has positive dimension, it has to contain non-CM points, whose Mumford--Tate group is not abelian.
By  Lemma~\ref{prop:tricotomy}, the  Lie subalgebra $\fg_{\Q_\ell}\cap \mathfrak{su}(3,\Q_\ell)$  contains a form of $\mathfrak{sl}(2,\Q_\ell)$. 
By  \cite[Thm.~1.1]{cadoret-tamagawa-duke1}  applied to the abelian family $f:A\to C$  there exist $B>0$ such that for all  $x\in C(k)$ outside a finite set  $C_\rho$ we have
  \begin{equation}\label{eq:uniform}
    [\Gamma_C: \Gamma_x]\leqslant B,
    \end{equation}
where  $\Gamma_x=\rho_{A_x,\ell}(\Gal_k)$  with $A_x$  the abelian $3$-fold of Picard type corresponding to $x$.

\begin{lemma}
There exists $r=r(C)\in \Z$ such that $\Gamma_x\supset \exp\left(\mathfrak{su}(2,\ell^{r}\Z_\ell)\right)$ for all $x\in C(k){\smallsetminus} C_\rho$.
\end{lemma}

\begin{proof}
We fix an exponential map on $\mathfrak{su}(2,\Q_\ell)$ so that  $\Gamma_C\supset \exp\left(\mathfrak{su}(2,\Z_\ell)\right)$. 
Using  that a subgroup of index at most $B$ contains a normal subgroup of index at most $B!$,  \eqref{eq:uniform} implies that 
 $\Gamma_x$ contains $B! \cdot \exp\left(\mathfrak{su}(2,\Z_\ell)\right)=\exp\left(\mathfrak{su}(2,\ell^{r}\Z_\ell)\right)$, where $r$ is the $\ell$-adic valuation of $B!$.
\end{proof}

It follows that for all $r\geqslant r(C)$ and for all $A$ as above  corresponding to a $k$-rational point on $C{\smallsetminus} C_\rho$,  $A[\ell^r]$ does not admit a  full $k$-rational flag. Finally, applying  Proposition~\ref{prop:MT} to the finitely many remaining  $k$-rational points,  yields  an integer $r$ such that all $k$-rational points in $Y_{\wK''_0(\ell^r)}$ are of CM type, completing the proof of  the Theorem.

\bigskip 
  { \noindent{\it Acknolwedgements:} { \small
We would like to thank our respective institutions, and also the TIFR where part of the work was done. We would like to thank A.~ Cadoret, T.~Graber, M.~ Goreski, B.~Gross, H.~ Hida, A.~Jorza, K.-W.~Lan, B.~Mazur, J.~ Nekov\'a\v{r}, D.~ Prasad, A.~ Raghuram, D.~Rohrlich, J.-P.~Serre and  J.~ Tilouine  for helpful discussions. }

\bibliographystyle{siam}

\end{document}